# EXTREMAL BEHAVIOR OF STOCHASTIC INTEGRALS DRIVEN BY REGULARLY VARYING LÉVY PROCESSES

By Henrik Hult[1] and Filip Lindskog

*Brown University and KTH Stockholm*

We study the extremal behavior of a stochastic integral driven by a multivariate Lévy process that is regularly varying with index $\alpha > 0$. For predictable integrands with a finite $(\alpha + \delta)$-moment, for some $\delta > 0$, we show that the extremal behavior of the stochastic integral is due to one big jump of the driving Lévy process and we determine its limit measure associated with regular variation on the space of càdlàg functions.

**1. Introduction.** Stochastic integrals driven by Lévy processes constitute a broad and popular class of semimartingales used as the driving noise in a wide variety of probabilistic models, for instance, the evolution of assets prices in mathematical finance. The extremal behavior of these processes is of importance when computing failure probabilities in various systems, for example, the probability that a functional of the sample path of the process exceeds some high threshold. In the presence of heavy tails of the underlying noise process such failures are often most likely due to one or a few unlikely events, such as large discontinuities (jumps) of the driving noise process. In the presence of Pareto-like tails of the underlying distributions regular variation on the space of càdlàg functions provides a useful framework to describe the extremal behavior of stochastic processes and approximate failure probabilities. In this paper we study the extremal behavior of stochastic integrals with respect to regularly varying Lévy processes. A first step toward studying the extremes of these processes was communicated to the authors by D. Applebaum [2].

The notion of regular variation is fundamental in various fields of applied probability. It serves as domain of attraction condition for partial sums of

Received March 2005; revised November 2005.
[1]Supported by the Swedish Research Council.

i.i.d. random vectors [28] or for componentwise maxima of i.i.d. random vectors [26], and it occurs in a natural way for the finite dimensional distributions of the stationary solution to stochastic recurrence equations [16, 22], including ARCH and GARCH processes; see [5], compare Section 8.4 in [14]. Let us consider an $\mathbb{R}^d$-valued vector $\mathbf{X}$. We call it *regularly varying* if there exists a sequence $(a_n)$ of positive numbers such that $a_n \uparrow \infty$ and a nonzero Radon measure $\mu$ on the $\sigma$-field $\mathcal{B}(\overline{\mathbb{R}}_0^d)$ of the Borel sets of $\overline{\mathbb{R}}_0^d = \overline{\mathbb{R}}^d \setminus \{\mathbf{0}\}$ (with $\overline{\mathbb{R}} = [-\infty, \infty]$) such that

$$(1.1) \qquad \mu(\overline{\mathbb{R}}^d \setminus \mathbb{R}^d) = 0 \quad \text{and} \quad n \, \mathrm{P}(a_n^{-1} \mathbf{X} \in \cdot) \xrightarrow{v} \mu(\cdot),$$

where $\xrightarrow{v}$ denotes vague convergence on $\mathcal{B}(\overline{\mathbb{R}}_0^d)$. We write $\mathbf{X} \in \mathrm{RV}((a_n), \mu, \overline{\mathbb{R}}_0^d)$. For details on the concept of vague convergence, we refer to [10], [21] and [26]. It can be shown that (1.1) necessarily implies that $\mu(uA) = u^{-\alpha} \mu(A)$ for some $\alpha > 0$, all $u > 0$ and all Borel sets $A$ bounded away from $\mathbf{0}$. Therefore, we also refer to *regular variation with index $\alpha$* in this context.

Definition (1.1) of regular variation has the advantage that it can be extended to random elements $\mathbf{X}$ with values in a separable Banach space [3] or certain linear metric spaces. We will use a formulation introduced in [11]. There the authors used regular variation of stochastic processes in the space of continuous functions and in the Skorokhod space $\mathbb{D}[0,1]$ in connection with max-stable distributions to extend many of the important results in classical extreme value theory to an infinite-dimensional setting. See also [15] for related results. This construction was taken up in [17], where regular variation of stochastic processes with values in the space $\mathbb{D} = \mathbb{D}([0,1], \mathbb{R}^d)$ of $\mathbb{R}^d$-valued càdlàg functions on $[0,1]$, equipped with the $J_1$-topology (see [6]), was considered. There regular variation of càdlàg processes was characterized in terms of regular variation of their finite dimensional distributions in the sense of (1.1) and a relative compactness condition in the spirit of weak convergence of stochastic processes [6]. Then an application of the continuous mapping theorem yields the tail behavior of interesting functionals.

In this paper we study the extremal behavior of a stochastic integral $(\mathbf{Y} \cdot \mathbf{X})$ given by

$$(1.2) \quad (\mathbf{Y} \cdot \mathbf{X})_t = \int_0^t \mathbf{Y}_s \, d\mathbf{X}_s = \left( \int_0^t Y_s^{(1)} \, dX_s^{(1)}, \ldots, \int_0^t Y_s^{(d)} \, dX_s^{(d)} \right),$$

$t \in [0,1]$, where $\mathbf{X} = (X_t^{(1)}, \ldots, X_t^{(d)})_{t \in [0,1]}$ is a $d$-dimensional Lévy process, which is regularly varying with index $\alpha > 0$. The stochastic process $\mathbf{Y} = (Y_t^{(1)}, \ldots, Y_t^{(d)})_{t \in [0,1]}$ is predictable càglàd and satisfies the moment condition $\mathrm{E}(\sup_{t \in [0,1]} |\mathbf{Y}_t|^{\alpha+\delta}) < \infty$ for some $\delta > 0$, where $|\cdot|$ denotes the Euclidean norm on $\mathbb{R}^d$.

It is known (see, e.g., [17]) that the extremal behavior of a multivariate regularly varying Lévy process is due to one large jump. Therefore, it is natural to guess that the extremal behavior of the stochastic integral (1.2) is due to one large jump of the underlying Lévy process. This is indeed the case. We begin by showing that (see Theorem 3.3), for each $\varepsilon > 0$ (with $|\mathbf{x}|_\infty = \sup_{t \in [0,1]} |\mathbf{x}_t|$),

$$\lim_{u \to \infty} P(d^\circ(u^{-1}\mathbf{X}, u^{-1}\Delta\mathbf{X}_\tau 1_{[\tau,1]}) > \varepsilon \mid |\mathbf{X}|_\infty > u) = 0,$$
(1.3)
$$\lim_{u \to \infty} P(d^\circ(u^{-1}\mathbf{X}, u^{-1}\Delta\mathbf{X}_\tau 1_{[\tau,1]}) > \varepsilon \mid |\Delta\mathbf{X}_\tau| > u) = 0,$$

where $d^\circ$ is the $J_1$-metric on the space of càdlàg functions, $\tau$ denotes the time of the jump of $\mathbf{X}$ with largest norm, and $\Delta\mathbf{X}_\tau = \mathbf{X}_\tau - \mathbf{X}_{\tau-}$. The interpretation of (1.3) is that when $\mathbf{X}$ is extreme (i.e., when $|\mathbf{X}|_\infty > u$ and $u$ is large) its sample path is well approximated (in an asymptotic sense) by a step function with one step. The second part of (1.3) implies that there is no other contribution to the extremal behavior of $\mathbf{X}$.

By the Lévy–Itô decomposition (e.g., [29], page 120), $\mathbf{X}$ can be decomposed into a sum of two independent processes

$$\mathbf{X} = \widetilde{\mathbf{X}} + \mathbf{J},\tag{1.4}$$

where $\mathbf{J}$ is a compound Poisson process with points $(\mathbf{Z}_k, \tau_k)$ and $|\mathbf{Z}_k| \geq 1$, that is,

$$\mathbf{J}_t = \sum_{k=1}^{N_t} \mathbf{Z}_k,$$

where $(N_t)$, given by $N_t = \sup\{k : \tau_k \leq t\}$, is a Poisson process. With this representation one can show that $\mathbf{X}$ is large, because one of the $\mathbf{Z}_k$'s is large whereas $\widetilde{\mathbf{X}}$ has light tails and does not have any influence on the extremal behavior of $\mathbf{X}$. Furthermore, the stochastic integral may be written as

$$(\mathbf{Y} \cdot \mathbf{X})_t = (\mathbf{Y} \cdot \widetilde{\mathbf{X}})_t + \sum_{k=1}^{N_t} \mathbf{Y}_{\tau_k} \mathbf{Z}_k.$$

Throughout the paper $\mathbf{xy}$ denotes componentwise multiplication, that is, $\mathbf{xy} = (x^{(1)}y^{(1)}, \ldots, x^{(d)}y^{(d)})$. If $\mathbf{Y}$ is predictable and $E(|\mathbf{Y}|_\infty^{\alpha+\delta}) < \infty$ for some $\delta > 0$, then it seems plausible, in the light of a classical result by Breiman [9] for the tail behavior of products of independent random variables, that $(\mathbf{Y} \cdot \mathbf{X})$ is well approximated by $\mathbf{Y}_\tau \mathbf{Z}_{k^*} 1_{[\tau,1]} = \mathbf{Y}_\tau \Delta\mathbf{X}_\tau 1_{[\tau,1]}$ given that $|(\mathbf{Y} \cdot \mathbf{X})|_\infty$ is large. Here $k^*$ denotes the index of the large jump, $\tau_{k^*} = \tau$. Indeed, Theorem 3.4 shows that

$$\lim_{u \to \infty} P(d^\circ(u^{-1}(\mathbf{Y} \cdot \mathbf{X}), u^{-1}\mathbf{Y}_\tau \Delta\mathbf{X}_\tau 1_{[\tau,1]}) > \varepsilon \mid |(\mathbf{Y} \cdot \mathbf{X})|_\infty > u) = 0,$$

$$\lim_{u \to \infty} P(d^\circ(u^{-1}(\mathbf{Y} \cdot \mathbf{X}), u^{-1}\mathbf{Y}_\tau \Delta\mathbf{X}_\tau 1_{[\tau,1]}) > \varepsilon \mid |\mathbf{Y}_\tau \Delta\mathbf{X}_\tau| > u) = 0.$$

Moreover, the process $(\mathbf{Y} \cdot \mathbf{X})$ is regularly varying on the space of càdlàg functions (see Section 2 for details). That is, there exist a limit measure $m^*$ and a sequence $(a_n)$ of positive numbers such that $a_n \uparrow \infty$ and for all sets $B \in \mathcal{B}(\mathbb{D})$ bounded away from $\mathbf{0}$ with $m^*(\partial B) = 0$, we have

$$n \operatorname{P}(a_n^{-1}(\mathbf{Y} \cdot \mathbf{X}) \in B) \to m^*(B).$$

We compute the limit measure $m^*$ as

$$m^*(B) = \operatorname{E}(\mu\{\mathbf{x} \in \overline{\mathbb{R}}_0^d : \mathbf{Y}_V \mathbf{x} 1_{[V,1]} \in B\}),$$

where $\mu$ is the regular variation limit measure (on $\overline{\mathbb{R}}_0^d$) of the jumps $\mathbf{Z}_k$ of the Lévy process and $V$ is uniformly distributed on $[0,1)$ and independent of the process $\mathbf{Y}$. As a simple illustration, we consider a univariate Lévy process $(X_t)_{t \in [0,1]}$ with $\operatorname{P}(X_1 > u) = u^{-\alpha} L(u)$ for some slowly varying function $L$. If $(Y_t)_{t \in [0,1]}$ is a non-negative process that satisfies the relevant conditions, then a straightforward application of Theorem 3.4 and the Continuous Mapping Theorem yields

$$\operatorname{P}\left(\int_0^t Y_v \, dX_v > u\right) \sim \operatorname{P}\left(\sup_{s \in [0,t]} \int_0^s Y_v \, dX_v > u\right)$$
$$\sim \operatorname{E}\left(\int_0^t Y_v^\alpha \, dv\right) u^{-\alpha} L(u),$$

where $f(u) \sim g(u)$ means that $\lim_{u \to \infty} f(u)/g(u) = 1$.

Stochastic integrals of the type (1.2) are encountered in many applications, in particular, in mathematical finance. Empirical evidence of regularly varying distributions in finance is recorded, for instance, in [1, 14] and [23]. In a financial context the process $\mathbf{Y}$ may be interpreted as a volatility process and the integral (1.2) the evolution of the log prices of $d$ assets.

The paper is organized as follows. In Section 2 we recall the concept of regular variation for stochastic processes with càdlàg sample paths (regular variation on $\mathbb{D}$). Section 3 contains the main results, which include an extension of Breiman's theorem to independent càdlàg processes, a result on approximating the trajectories of regularly varying processes, and the main theorem of this paper concerning the extremal behavior of stochastic integrals. The remaining Sections 4 and 5 contain the proofs and some auxiliary results.

**2. Regular variation and Lévy processes.** Let us recall the notion of regular variation for stochastic processes with sample paths in $\mathbb{D} = \mathbb{D}([0,1], \mathbb{R}^d)$; the space of functions $\mathbf{x} : [0,1] \to \mathbb{R}^d$ that are right continuous with left limits. This space is equipped with the so-called $J_1$-metric (referred to as $d^\circ$ in [6]) that makes it complete and separable.

We denote by $\mathbb{S}_\mathbb{D}$ the subspace $\{\mathbf{x} \in \mathbb{D} : |\mathbf{x}|_\infty = 1\}$ (where $|\mathbf{x}|_\infty = \sup_{t \in [0,1]} |\mathbf{x}_t|$) equipped with the subspace topology. Define $\overline{\mathbb{D}}_0 = (0, \infty] \times \mathbb{S}_\mathbb{D}$, where $(0, \infty]$ is equipped with the metric $\rho(x, y) = |1/x - 1/y|$, making it complete and separable. Then the space $\overline{\mathbb{D}}_0$, equipped with the metric $\max\{\rho(x^*, y^*), d^\circ(\widetilde{\mathbf{x}}, \widetilde{\mathbf{y}})\}$, is a complete separable metric space. For $\mathbf{x} = (x^*, \widetilde{\mathbf{x}}) \in \overline{\mathbb{D}}_0$, we write $|\mathbf{x}|_\infty = x^*$. The topological spaces $\mathbb{D}\backslash\{\mathbf{0}\}$ (equipped with the subspace topology of $\mathbb{D}$) and $(0, \infty) \times \mathbb{S}_\mathbb{D}$ (equipped with the subspace topology of $\overline{\mathbb{D}}_0$) are homeomorphic; the mapping $T$ given by $T(\mathbf{x}) = (|\mathbf{x}|_\infty, \mathbf{x}/|\mathbf{x}|_\infty)$ is a homeomorphism. Hence,

$$\mathcal{B}(\overline{\mathbb{D}}_0) \cap ((0, \infty) \times \mathbb{S}_\mathbb{D}) = \mathcal{B}(T(\mathbb{D}\backslash\{\mathbf{0}\})),$$

that is, the Borel sets of $\mathcal{B}(\overline{\mathbb{D}}_0)$ that are of interest to us can be identified with the usual Borel sets on $\mathbb{D}$ (viewed in polar coordinates) that do not contain the zero function. For notational convenience, we will throughout the paper identify $\mathbb{D}$ with the product space $[0, \infty) \times \mathbb{S}_\mathbb{D}$ so that expressions like $\overline{\mathbb{D}}_0 \backslash \mathbb{D}$ $(= \{\infty\} \times \mathbb{S}_\mathbb{D})$ make sense.

Regular variation on $\mathbb{D}$ is naturally expressed in terms of so-called $\hat{w}$-convergence of boundedly finite measures on $\overline{\mathbb{D}}_0$. A boundedly finite measure assigns finite measure to bounded sets. A sequence of boundedly finite measures $(m_n)_{n \in \mathbb{N}}$ on a complete separable metric space $\mathbb{E}$ converges to $m$ in the $\hat{w}$-topology, $m_n \xrightarrow{\hat{w}} m$, if $m_n(B) \to m(B)$ for every bounded Borel set $B$ with $m(\partial B) = 0$. If the state space $\mathbb{E}$ is locally compact, which $\overline{\mathbb{D}}_0$ is not but $\overline{\mathbb{R}}_0^d$ ($\overline{\mathbb{R}} = [-\infty, \infty]$) is, then a boundedly finite measure is called a Radon measure, and $\hat{w}$-convergence coincides with vague convergence and we write $m_n \xrightarrow{v} m$. Finally we note that if $m_n \xrightarrow{\hat{w}} m$ and $m_n(\mathbb{E}) \to m(\mathbb{E}) < \infty$, then $m_n \xrightarrow{w} m$ with $\xrightarrow{w}$ denoting weak convergence. For details on $\hat{w}$-, vague- and weak convergence, we refer to [10], Appendix 2.

Recall the definition (1.1) of multivariate regular variation. For a stochastic process with sample paths in $\mathbb{D}$, regular variation can be formulated similarly. A stochastic process $\mathbf{X} = (\mathbf{X}_t)_{t \in [0,1]}$ with sample paths in $\mathbb{D}$ is said to be regularly varying if there exist a sequence $(a_n)$, $0 < a_n \uparrow \infty$, and a nonzero boundedly finite measure $m$ on $\mathcal{B}(\overline{\mathbb{D}}_0)$ with $m(\overline{\mathbb{D}}_0 \backslash \mathbb{D}) = 0$ such that, as $n \to \infty$,

$$n \mathrm{P}(a_n^{-1} \mathbf{X} \in \cdot) \xrightarrow{\hat{w}} m(\cdot) \qquad \text{on } \mathcal{B}(\overline{\mathbb{D}}_0).$$

We write $\mathbf{X} \in \mathrm{RV}((a_n), m, \overline{\mathbb{D}}_0)$. If $\nu$ is a measure satisfying, with $(a_n)$ and $m$ as above, $n\nu(a_n \cdot) \xrightarrow{\hat{w}} m(\cdot)$ on $\mathcal{B}(\overline{\mathbb{D}}_0)$, then we write $\nu \in \mathrm{RV}((a_n), m, \overline{\mathbb{D}}_0)$ and similarly for measures on $\mathbb{R}^d$.

REMARK 2.1. (i) Theorem 10 in [17] gives necessary and sufficient conditions for $\mathbf{X} \in \mathrm{RV}((a_n), m, \overline{\mathbb{D}}_0)$ in terms of multivariate regular variation for finite dimensional distributions of $\mathbf{X}$ and a relative compactness condition.

(ii) If $\mathbf{X} \in \mathrm{RV}((a_n), m, \overline{\mathbb{D}}_0)$, then there exists $\alpha > 0$ such that $m(u\cdot) = u^{-\alpha} m(\cdot)$ for every $u > 0$ (e.g., [19], Theorem 3.1). Therefore, we will also refer to regular variation with index $\alpha > 0$ or $\mathbf{X} \in \mathrm{RV}_\alpha((a_n), m, \overline{\mathbb{D}}_0)$.

For other equivalent formulations of regular variation on $\overline{\mathbb{R}}_0^d$ (most of which can be modified into formulations of regular variation on $\overline{\mathbb{D}}_0$), we refer to [4, 5, 19, 26, 27]. For the classical theory of regularly varying functions, see [8].

The next theorem is an analogue of the Continuous Mapping Theorem for weak convergence. Let $\mathrm{Disc}(h)$ denote the set of discontinuities of a mapping $h$ from a metric space $\mathbb{E}$ to a metric space $\mathbb{E}'$. It is shown on page 225 in [6] that $\mathrm{Disc}(h) \in \mathcal{B}(\mathbb{E})$.

THEOREM 2.1. *Let $\mathbf{X} = (\mathbf{X}_t)_{t \in [0,1]}$ be a stochastic process with sample paths in $\mathbb{D}$ and let $\mathbb{E}'$ be a complete separable metric space. Suppose that $\mathbf{X} \in \mathrm{RV}((a_n), m, \overline{\mathbb{D}}_0)$ and that $h : \overline{\mathbb{D}}_0 \to \mathbb{E}'$ is a measurable mapping satisfying $m(\mathrm{Disc}(h)) = 0$ and $h^{-1}(B)$ is bounded in $\overline{\mathbb{D}}_0$ for every bounded $B \in \mathcal{B}(\mathbb{E}')$. Then, as $n \to \infty$,*

$$n \, \mathrm{P}(h(a_n^{-1} \mathbf{X}) \in \cdot) \stackrel{\hat{w}}{\to} m \circ h^{-1}(\cdot) \qquad \text{on } \mathcal{B}(\mathbb{E}').$$

See [17], Theorem 6, for a proof.

REMARK 2.2. The conclusion of the theorem holds for random vectors, that is, if one considers $\mathbf{X} \in \mathrm{RV}((a_n), m, \overline{\mathbb{R}}_0^d)$ and mappings $h : \overline{\mathbb{R}}_0^d \to \mathbb{E}'$.

Given a regularly varying stochastic process $\mathbf{X}$ with limit measure $m$, the continuous mapping theorem allows us to derive the asymptotic behavior of mappings $h(\mathbf{X})$ of the sample paths, for instance, the componentwise supremum and average;

$$\left( \sup_{t \in [0,1]} |X_t^{(1)}|, \ldots, \sup_{t \in [0,1]} |X_t^{(d)}| \right) \quad \text{and} \quad \left( \int_0^1 X_s^{(1)} \, ds, \ldots, \int_0^1 X_s^{(d)} \, ds \right).$$

Thus, if we are interested in approximating the failure probability of a certain regularly varying stochastic process $\mathbf{X}$, expressed as the probability that $h(\mathbf{X})$ is in some set far away from the origin, then a natural approach is to first determine the limit measure $m$ of the processes and then apply the continuous mapping theorem. This is the reason for our interest in finding the limit measure for various regularly varying stochastic processes.

In the rest of this paper we will focus on the computation of the limit measure of a stochastic integral with respect to a (multivariate) Lévy process. We first recall some relevant results on regular variation of a Lévy process and, more generally, of Markov processes with increments satisfying a condition of weak dependence (see [17]).

We will frequently use the Lévy–Itô decomposition (e.g., [29], page 120) which says that a Lévy process $\mathbf{X}$ on $\mathbb{R}^d$ with generating triplet $(A, \boldsymbol{\gamma}, \nu)$ may be decomposed as

$$\mathbf{X} = \widetilde{\mathbf{X}} + \mathbf{J} \quad \text{a.s.,} \tag{2.1}$$

where, for almost all $\omega \in \Omega$,

$$\widetilde{\mathbf{X}}_t(\omega) = \lim_{\gamma \to 0} \int_{(0,t] \times \{\gamma \leq |\mathbf{x}| < 1\}} \mathbf{x}\{\xi(d(s,\mathbf{x}),\omega) - ds\nu(d\mathbf{x})\} + \boldsymbol{\gamma} t + \mathbf{W}_t(\omega), \tag{2.2}$$

$$\mathbf{J}_t(\omega) = \int_{(0,t] \times \{|\mathbf{x}| \geq 1\}} \mathbf{x}\xi(d(s,\mathbf{x}),\omega), \tag{2.3}$$

$\xi$ is a Poisson random measure with mean measure $\lambda \times \nu$ [$\xi \sim \text{PRM}(\lambda \times \nu)$, $\lambda$ denoting Lebesgue measure], and $\mathbf{W}$ is a Gaussian process with stationary and independent increments. The processes $\widetilde{\mathbf{X}}$ and $\mathbf{J}$ are independent.

For a Lévy process $\mathbf{X}$ regular variation on $\mathbb{D}$ is intimately connected to regular variation of the Lévy measure $\nu$ of $\mathbf{X}_1$. This is summarized in the following lemma.

LEMMA 2.1. *Let $\mathbf{X}$ be a Lévy process with Lévy measure $\nu$. The following statements are equivalent:*

  (i) $\mathbf{X}_1 \in \text{RV}((a_n), \mu, \overline{\mathbb{R}}_0^d)$,
  (ii) $\nu \in \text{RV}((a_n), \mu, \overline{\mathbb{R}}_0^d)$,
  (iii) $\mathbf{X} \in \text{RV}((a_n), m, \overline{\mathbb{D}}_0)$ *with* $m_t = t\mu$ *for every* $t \in [0,1]$.

The proof of these statements follows by combining Proposition 3.1 in [18] and Theorem 10 in [17]. In the univariate case ($d = 1$) a proof of the equivalence (i) $\Leftrightarrow$ (ii) was given in [13]. The limit measure $m$ in (iii) is concentrated on the set of step functions with one step; that is, $m(\mathcal{V}^c) = 0$, where $\mathcal{V}^c$ is the complement of

$$\mathcal{V} = \{\mathbf{x} \in \mathbb{D} : \mathbf{x} = \mathbf{z}1_{[v,1]}, v \in [0,1), \mathbf{z} \in \mathbb{R}^d \setminus \{\mathbf{0}\}\} \tag{2.4}$$

(see [17], Theorem 15). Moreover, the measure $m$ has the representation (see [20], Remark 2.1)

$$m(B) = \int_{[0,1]} \int_{\mathbb{R}_0^d} 1_B(\mathbf{y}1_{[t,1]}) \mu(d\mathbf{y}) \, dt, \tag{2.5}$$

where $\mathbf{y}1_{[t,1]}$ is the element $f \in \mathbb{D}$ given by $f(u) = 0$ for $u \in [0,t)$ and $f(u) = \mathbf{y}$ for $u \in [t,1]$.

**3. Main results.** We assume that all random elements are defined on a filtered complete probability space $(\Omega, \mathcal{F}, (\mathcal{F}_t)_{t \in [0,1]}, P)$ satisfying the usual hypotheses (see [25], page 3).

3.1. *Regular variation for products of independent stochastic processes.* Before we study the stochastic integral in more detail in Section 3.3, we first consider a much simpler situation; products of independent stochastic processes. In this section we will extend a well-known result by Breiman [9], Proposition 3, concerning the tail behavior of products of independent random variables to stochastic processes with sample paths in $\mathbb{D}$. Breiman's result (more precisely, a slight generalization of this result) says that for independent nonnegative random variables $Y$ and $X$ such that $X$ is regularly varying with index $\alpha$ and $E(Y^{\alpha+\delta}) < \infty$ for some $\delta > 0$, as $x \to \infty$,

$$[P(X > x)]^{-1} P(YX > x) \to E(Y^\alpha).$$

Since regular variation of $X$ can be formulated in terms of vague convergence on $(0, \infty]$, there exist a sequence $(a_n)$, $0 < a_n \uparrow \infty$, and a nonzero Radon measure $\mu$ on $\mathcal{B}((0, \infty])$ such that, as $n \to \infty$,

$$n P(a_n^{-1} X \in \cdot) \xrightarrow{v} \mu(\cdot) \qquad \text{on } \mathcal{B}((0, \infty]),$$

and $\mu((u, \infty]) = cu^{-\alpha}$. Then Breiman's result may be written as

$$
\begin{aligned}
n P(a_n^{-1} Y X \in \cdot) &\xrightarrow{v} E(\mu\{x \in (0, \infty] : Yx \in \cdot\}) \\
&= E(Y^\alpha) \mu(\cdot) \qquad \text{on } \mathcal{B}((0, \infty]).
\end{aligned}
$$ (3.1)

This result was extended to regularly varying random vectors in [5], Proposition A.1. Our version of Breiman's result for stochastic processes is Theorem 3.1 below. Given an element $\mathbf{y} \in \mathbb{D}$, let $\phi_\mathbf{y} : \mathbb{D} \to \mathbb{D}$ be given by

$$\phi_\mathbf{y}(\mathbf{x}) = \mathbf{y}\mathbf{x} = (y^{(1)} x^{(1)}, \ldots, y^{(d)} x^{(d)}). \tag{3.2}$$

Then $\phi_\mathbf{y}$ is measurable and continuous at those $\mathbf{x}$ for which $\text{Disc}(\mathbf{x}) \cap \text{Disc}(\mathbf{y}) = \varnothing$ (see [30]).

THEOREM 3.1. *Let $\mathbf{X}$ and $\mathbf{Y}$ be independent stochastic processes with sample paths in $\mathbb{D}$. Suppose that $\mathbf{X} \in \text{RV}_\alpha((a_n), m, \overline{\mathbb{D}}_0)$, that $E(|\mathbf{Y}|_\infty^{\alpha+\delta}) < \infty$ for some $\delta > 0$ and that $\min_{k=1,\ldots,d} |Y^{(k)}|_\infty > 0$ a.s. If $E(m(\text{Disc}(\phi_\mathbf{Y}))) = 0$, then, as $n \to \infty$,*

$$n P(a_n^{-1} \mathbf{Y} \mathbf{X} \in \cdot) \xrightarrow{\hat{w}} E(m \circ \phi_\mathbf{Y}^{-1}(\cdot)) = E(m\{\mathbf{x} \in \overline{\mathbb{D}}_0 : \mathbf{Y}\mathbf{x} \in \cdot\}) \qquad \text{on } \mathcal{B}(\overline{\mathbb{D}}_0).$$

REMARK 3.1. (i) If the process $\mathbf{Y}$ has continuous sample paths a.s., then $E(m(\text{Disc}(\phi_\mathbf{Y}))) = 0$.

(ii) If $\mathbf{X}$ is a Lévy process, then $\mathrm{E}(m(\mathrm{Disc}(\phi_{\mathbf{Y}}))) = 0$ for all càdlàg processes $\mathbf{Y}$ (see Lemma 5.1).

(iii) If $\mathrm{E}(m \circ \phi_{\mathbf{Y}}^{-1}(\cdot))$ is a nonzero measure, then $\mathbf{YX}$ is regularly varying, that is, $\mathbf{YX} \in \mathrm{RV}_\alpha((a_n), \mathrm{E}(m \circ \phi_{\mathbf{Y}}^{-1}(\cdot)), \overline{\mathbb{D}}_0)$.

(iv) If $\mathbf{X}$ is a Lévy process, then, by (2.5),

$$\mathrm{E}(m \circ \phi_{\mathbf{Y}}^{-1}(\{\mathbf{z} \in \overline{\mathbb{D}}_0 : |\mathbf{z}|_\infty \geq 1\})) = \mathrm{E}(\mu\{\mathbf{x} \in \overline{\mathbb{R}}_0^d : |\mathbf{Yx}1_{[V,1]}|_\infty \geq 1\}),$$

where $V$ is uniformly distributed on $[0, 1)$ and independent of $\mathbf{Y}$. Set

$$A_{\varepsilon,\delta} = \left\{\omega \in \Omega : \min_{k=1,\ldots,d} \sup_{t \in [0, 1-\delta]} |Y_t^{(k)}(\omega)| > \varepsilon\right\}.$$

By assumption and since $\mathbf{Y}$ has right-continuous sample paths, there exist $\delta \in (0, 1)$ and $\varepsilon, \eta > 0$ such that $\mathrm{P}(A_{\varepsilon,\delta}) > \eta$. Hence,

$$\mathrm{E}(\mu\{\mathbf{x} \in \overline{\mathbb{R}}_0^d : |\mathbf{Yx}1_{[V,1]}|_\infty \geq 1\})$$
$$\geq \mathrm{E}(\mu\{\mathbf{x} \in \overline{\mathbb{R}}_0^d : |\mathbf{Yx}1_{[V,1]}|_\infty \geq 1\}; A_{\varepsilon,\delta}, V \leq \delta)$$
$$> \mu\left\{\mathbf{x} \in \overline{\mathbb{R}}_0^d : \max_{k=1,\ldots,d} |x^{(k)}| > 1/\varepsilon\right\} \eta \delta > 0.$$

Hence, $\mathrm{E}(m \circ \phi_{\mathbf{Y}}^{-1}(\cdot))$ is a nonzero measure for every $\mathbf{Y}$ satisfying the conditions of the theorem.

REMARK 3.2. It was shown by Embrechts and Goldie [12], corollary to Theorem 3, that for nonnegative random variables $X$ and $Y$ with $X$ regularly varying with index $\alpha$ and $\mathrm{P}(Y > x) = o(\mathrm{P}(X > x))$ as $x \to \infty$, it holds that $YX$ is regularly varying with index $\alpha$. However, for independent random vectors $\mathbf{X}$ and $\mathbf{Y}$, regular variation of $\mathbf{X}$ and $\mathrm{P}(|\mathbf{Y}| > x) = o(\mathrm{P}(|\mathbf{X}| > x))$ as $x \to \infty$ is not sufficient for regular variation of $\mathbf{YX}$. Therefore, it is, in general, not possible to replace the moment condition in Theorem 3.1 by $\mathrm{P}(|\mathbf{Y}|_\infty > x) = o(\mathrm{P}(|\mathbf{X}|_\infty > x))$ as $x \to \infty$.

3.2. *Approximating the extreme sample paths of regularly varying stochastic processes.* As explained in Section 2, the limit measure associated with regular variation of a stochastic process in $\mathbb{D}$ characterizes its extremal behavior. Moreover, the continuous mapping theorem can be applied to derive the tail behavior of functionals of its sample paths. However, these results concern only the distributional aspects of the extremal behavior. In some cases we would like stronger results on approximating the extremal behavior of a stochastic process. We take the following approach. Consider two stochastic processes $\mathbf{X}$ and $\mathbf{Y}$ with sample paths in $\mathbb{D}$. If, given that $\mathbf{Y}$ is extreme (i.e., $|\mathbf{Y}|_\infty > u$ for $u$ large), the distance between the rescaled processes $u^{-1}\mathbf{X}$ and $u^{-1}\mathbf{Y}$ is small with high probability, then the extreme sample path behavior of $\mathbf{Y}$ may be approximated by that of $\mathbf{X}$. To conclude that there is no other contribution to the extreme sample paths of $\mathbf{Y}$ we also need that the distance between $u^{-1}\mathbf{X}$ and $u^{-1}\mathbf{Y}$ is small when $|\mathbf{X}|_\infty > u$ for large $u$. We say that the extreme sample paths of $\mathbf{Y}$ can be approximated by those of $\mathbf{X}$ and vice versa if, for every $\varepsilon > 0$,

$$
\begin{aligned}
\lim_{u\to\infty} \mathrm{P}(d^\circ(u^{-1}\mathbf{X}, u^{-1}\mathbf{Y}) > \varepsilon \mid |\mathbf{Y}|_\infty > u) = 0, \\
\lim_{u\to\infty} \mathrm{P}(d^\circ(u^{-1}\mathbf{X}, u^{-1}\mathbf{Y}) > \varepsilon \mid |\mathbf{X}|_\infty > u) = 0.
\end{aligned}
\tag{3.3}
$$

We typically look for a simple process $\mathbf{X}$ (e.g., a step function) such that (3.3) holds.

Theorem 3.2 below says that if (3.3) holds and $\mathbf{X}$ is regularly varying, then $\mathbf{Y}$ is regularly varying with the same limit measure. It is similar in spirit to the following well-known result for weak convergence: If $(\mathbb{E}, \rho)$ is a metric space and $(X_n, Y_n)$ are random elements of $\mathbb{E} \times \mathbb{E}$, then $X_n \xrightarrow{d} X$ and $\rho(X_n, Y_n) \xrightarrow{d} 0$ imply $Y_n \xrightarrow{d} X$ (see, e.g., [7], Theorem 3.1).

THEOREM 3.2. *Let $\mathbf{X}$ and $\mathbf{Y}$ be stochastic processes with sample paths in $\mathbb{D}$. If $\mathbf{X} \in \mathrm{RV}((a_n), m, \overline{\mathbb{D}}_0)$ and (3.3) holds, then $\mathbf{Y} \in \mathrm{RV}((a_n), m, \overline{\mathbb{D}}_0)$.*

Next we consider a regularly varying Lévy process $\mathbf{X} \in \mathrm{RV}((a_n), m, \overline{\mathbb{D}}_0)$. Let $\mathcal{V} \subset \mathbb{D}$ be the family of step functions in $\mathbb{D}$ with one step in (2.4). As already mentioned, the limit measure $m$ puts all its mass on this set. The next theorem is a slightly stronger version of this result: it describes, in the sense of (3.3), the sample paths of $\mathbf{X}$ given that $|\mathbf{X}|_\infty > u$ for $u$ large. First we need some notation. Define $\tau : \mathbb{D} \to [0, 1]$ as the time of the jump with largest norm of an element $\mathbf{x} \in \mathbb{D}$. If there are several jumps of equal size, we let $\tau(\mathbf{x})$ denote the first of them. More precisely,

$$
(3.4)\quad \tau(\mathbf{x}) = \lim_{\varepsilon \downarrow 0} \inf\{t \in (0,1) : |\Delta \mathbf{x}_t| = \sup\{|\Delta \mathbf{x}_s| : s \in (0,1), |\Delta \mathbf{x}_s| > \varepsilon\}\}.
$$

If the set in (3.4) is empty, then we put $\tau(\mathbf{x}) = 1$. The next result says that Lévy process $\mathbf{X} \in \mathrm{RV}((a_n), m, \overline{\mathbb{D}}_0)$ is asymptotically close to the step function given by $\Delta \mathbf{X}_{\tau(\mathbf{X})} 1_{[\tau(\mathbf{X}), 1]}$ in the sense of (3.3).

THEOREM 3.3. *Let $\mathbf{X} \in \mathrm{RV}((a_n), m, \overline{\mathbb{D}}_0)$ be a Lévy process. Then, for every $\varepsilon > 0$,*

$$
\begin{aligned}
\lim_{u\to\infty} \mathrm{P}(d^\circ(u^{-1}\mathbf{X}, u^{-1}\Delta\mathbf{X} 1_{[\tau(\mathbf{X}),1]}) > \varepsilon \mid |\mathbf{X}|_\infty > u) = 0, \\
\lim_{u\to\infty} \mathrm{P}(d^\circ(u^{-1}\mathbf{X}, u^{-1}\Delta\mathbf{X}_{\tau(\mathbf{X})} 1_{[\tau(\mathbf{X}),1]}) > \varepsilon \mid |\Delta\mathbf{X}_{\tau(\mathbf{X})}| > u) = 0
\end{aligned}
\tag{3.5}
$$

*and $\Delta\mathbf{X} 1_{[\tau(\mathbf{X}),1]} \in \mathrm{RV}((a_n), m, \overline{\mathbb{D}}_0)$.*

REMARK 3.3. If $\mathbf{X}$ and $\mathbf{Y}$ satisfy the hypotheses of Theorem 3.1, then one can also show that

$$\lim_{u\to\infty} P(d^\circ(u^{-1}\mathbf{YX}, u^{-1}\mathbf{Y}\Delta\mathbf{X}_\tau 1_{[\tau,1]}) > \varepsilon \mid |\mathbf{YX}|_\infty > u) = 0,$$

$$\lim_{u\to\infty} P(d^\circ(u^{-1}\mathbf{YX}, u^{-1}\mathbf{Y}\Delta\mathbf{X}_\tau 1_{[\tau,1]}) > \varepsilon \mid |\mathbf{Y}\Delta\mathbf{X}_\tau 1_{[\tau,1]}|_\infty > u) = 0,$$

with $\tau = \tau(\mathbf{X})$.

3.3. *Extremal behavior of stochastic integrals.* The main result in this paper concerns the extremal behavior of a stochastic integral $(\mathbf{Y} \cdot \mathbf{X})$ given by

$$(\mathbf{Y} \cdot \mathbf{X})_t = \left(\int_0^t Y_s^{(1)} \, dX_s^{(1)}, \ldots, \int_0^t Y_s^{(d)} \, dX_s^{(d)}\right), \qquad t \in [0,1],$$

where $\mathbf{X} \in \mathrm{RV}_\alpha((a_n), m, \overline{\mathbb{D}}_0)$ is a regularly varying Lévy process and $\mathbf{Y}$ is an $\mathbb{R}^d$-valued predictable càglàd process that satisfies the moment condition $\mathrm{E}(|\mathbf{Y}|_\infty^{\alpha+\delta}) < \infty$, for some $\delta > 0$. We refer to [25] for an account on stochastic integration. The intuitive idea is the following. Given that $|\mathbf{X}|_\infty$ is large, Theorem 3.3 states that $\mathbf{X}$ and $\Delta\mathbf{X}_\tau 1_{[\tau,1]}$ are asymptotically close, that is,

$$\mathbf{X} \approx \Delta\mathbf{X}_\tau 1_{[\tau,1]},$$

where $\tau = \tau(\mathbf{X})$ is the time of the jump with largest norm. This suggests that, given that $|(\mathbf{Y} \cdot \mathbf{X})|_\infty$ is large, we can replace $\mathbf{X}$ by $\Delta\mathbf{X}_\tau 1_{[\tau,1]}$ in the stochastic integral and thereby justify the following approximation, in the sense of (3.3):

$$(\mathbf{Y} \cdot \mathbf{X}) \approx \mathbf{Y}_\tau \Delta\mathbf{X}_\tau 1_{[\tau,1]}.$$

We have the following result:

THEOREM 3.4. *Let $\mathbf{X}$ be a Lévy process satisfying $\mathbf{X}_1 \in \mathrm{RV}_\alpha((a_n), \mu, \overline{\mathbb{R}}_0^d)$ and let $\mathbf{Y}$ be a predictable càglàd process satisfying $\mathrm{E}(|\mathbf{Y}|_\infty^{\alpha+\delta}) < \infty$ for some $\delta > 0$ and $\min_{k=1,\ldots,d} |Y^{(k)}|_\infty > 0$ a.s. Then, for every $\varepsilon > 0$,*

(3.6) $\quad \lim_{u\to\infty} P(d^\circ(u^{-1}(\mathbf{Y} \cdot \mathbf{X}), u^{-1}\mathbf{Y}_\tau \Delta\mathbf{X}_\tau 1_{[\tau,1]}) > \varepsilon \mid |(\mathbf{Y} \cdot \mathbf{X})|_\infty > u) = 0,$

(3.7) $\quad \lim_{u\to\infty} P(d^\circ(u^{-1}(\mathbf{Y} \cdot \mathbf{X}), u^{-1}\mathbf{Y}_\tau \Delta\mathbf{X}_\tau 1_{[\tau,1]}) > \varepsilon \mid |\mathbf{Y}_\tau \Delta\mathbf{X}_\tau| > u) = 0,$

*where $\tau = \tau(\mathbf{X})$. Moreover, $(\mathbf{Y} \cdot \mathbf{X}), \mathbf{Y}_\tau \Delta\mathbf{X}_\tau 1_{[\tau,1]} \in \mathrm{RV}_\alpha((a_n), m^*, \overline{\mathbb{D}}_0)$ with*

$$m^*(B) = \mathrm{E}(\mu\{\mathbf{x} \in \overline{\mathbb{R}}_0^d : \mathbf{Y}_V \mathbf{x} 1_{[V,1]} \in B\}),$$

*where $V$ is uniformly distributed on $[0,1)$ and independent of $\mathbf{Y}$.*

The idea behind the proof is the following. Using the Lévy–Itô decomposition (2.1), we can write $(\mathbf{Y} \cdot \mathbf{X}) = (\mathbf{Y} \cdot \mathbf{J}) + (\mathbf{Y} \cdot \widetilde{\mathbf{X}})$. Using the fact that $\widetilde{\mathbf{X}}$ has finite moments of all orders, we find that the extremal behavior will be determined by that of

$$(\mathbf{Y} \cdot \mathbf{J})_t = \sum_{k=1}^{N_t} \mathbf{Y}_{\tau_k} \mathbf{Z}_k,$$

where $(\mathbf{Z}_k)$ is an i.i.d. sequence with $\mathbf{Z}_k \in \mathrm{RV}((a_n), \mu, \overline{\mathbb{R}}_0^d)$ and independent of the Poisson process $(N_t)$. Since $\mathbf{Y}$ is predictable and $\tau_k$ is a stopping time, $\mathbf{Y}_{\tau_k}$ and $\mathbf{Z}_k$ are independent. Because of the moment condition, the multivariate version of Breiman's result gives the tail behavior of the product $\mathbf{Y}_{\tau_k}\mathbf{Z}_k$. Moreover, since the $\mathbf{Z}_k$'s are i.i.d. and $(N_t)$ is a Poisson process, we expect that asymptotically only one of the $\mathbf{Z}_k$'s will be large and, hence, that one term $\mathbf{Y}_{\tau_k}\mathbf{Z}_k$ will dominate the sum of the rest, that is, the extremal behavior of $(\mathbf{Y} \cdot \mathbf{J})$ is determined by $\mathbf{Y}_{\tau_{k^*}}\mathbf{Z}_{k^*}$, where $k^*$ is the index of the $\mathbf{Z}_k$'s with largest norm. The main difficulty comes from the fact that the terms $\mathbf{Y}_{\tau_k}\mathbf{Z}_k$ may be dependent. Note that since we only require that $\mathbf{Y}$ is predictable, $\mathbf{Y}_{\tau_k}$ may depend on the variables $\tau_1, \ldots, \tau_{k-1}$ and $\mathbf{Z}_1, \ldots, \mathbf{Z}_{k-1}$, as well as on $(\mathbf{Y}_s; s < \tau_k)$. To overcome this difficulty, we need a number of technical lemmas presented in Section 5. The limit measure for the stochastic integral $(\mathbf{Y} \cdot \mathbf{J})$ is computed in Proposition 5.1.

Let us now consider a couple of simple univariate examples that illustrate some of the applications of Theorem 3.4.

EXAMPLE 3.1. Let $X$ be a Lévy process with $X_1 \in \mathrm{RV}_\alpha((a_n), \mu, \overline{\mathbb{R}}_0)$ and with $\mu((u, \infty)) = cu^{-\alpha}$ for some $c > 0$. Let $Y$ satisfy the conditions of Theorem 3.4. If $Y_t > 0$ for all $t$, we may think of $Y$ as a volatility process and $(Y \cdot X)_t$ as the logarithm of an asset price at time $t$. Then $(Y \cdot X) \in \mathrm{RV}_\alpha((a_n), m^*, \overline{\mathbb{D}}_0)$, where $m^*$ is given by

$$m^*(B) = \mathrm{E}(\mu\{x \in \overline{\mathbb{R}}_0 : xY_V 1_{[V,1]} \in B\}),$$

and $V$ is uniformly distributed on $[0, 1)$ and independent of $Y$. In particular, applying the continuous mapping theorem with the functional $\pi_t : \mathbb{D} \to \mathbb{R}$ given by $\pi_t(z) = z_t$, we obtain, for each $u > 0$,

$$\begin{aligned}
n\,\mathrm{P}(a_n^{-1}(Y \cdot X)_t > u) &= n\,\mathrm{P}(a_n^{-1}(Y \cdot X) \in \pi_t^{-1}((u, \infty))) \\
&\to \mathrm{E}(\mu\{x \in \overline{\mathbb{R}}_0 : Y_V x 1_{[V,1]} \in \pi_t^{-1}((u, \infty))\}) \\
&= \mathrm{E}(\mu\{x \in \overline{\mathbb{R}}_0 : Y_V x > u\} 1_{[0,t]}(V)) \\
&= \mathrm{E}(Y_V^\alpha 1_{[0,t]}(V))\mu((u, \infty)) \\
&= c \int_0^t \mathrm{E}(Y_s^\alpha)\,ds\, u^{-\alpha}.
\end{aligned}$$

EXAMPLE 3.2. Consider the previous example and the supremum-functional $h_t : \mathbb{D} \to \mathbb{R}$ given by $h_t(z) = \sup_{s \in [0,t]} z_s$. We obtain, for each $u > 0$,

$$\begin{aligned}
n\,\mathrm{P}\!\left(a_n^{-1} \sup_{s \in [0,t]}(Y \cdot X)_s > u\right) &= n\,\mathrm{P}(a_n^{-1}(Y \cdot X) \in h_t^{-1}((u, \infty))) \\
&\to \mathrm{E}(\mu\{x \in \overline{\mathbb{R}}_0 : Y_V x 1_{[V,1]} \in h_t^{-1}((u, \infty))\}) \\
&= \mathrm{E}(\mu\{x \in \overline{\mathbb{R}}_0 : Y_V x > u\} 1_{[0,t]}(V)) \\
&= \mathrm{E}(Y_V^\alpha 1_{[0,t]}(V)) \mu((u, \infty)) \\
&= c \int_0^t \mathrm{E}(Y_s^\alpha)\, ds\, u^{-\alpha}.
\end{aligned}$$

As a consequence, we obtain that

$$\lim_{u \to \infty} [\mathrm{P}((Y \cdot X)_t > u)]^{-1} \mathrm{P}\!\left(\sup_{s \in [0,t]} (Y \cdot X)_s > u\right) = 1.$$

This extends the tail-equivalence for heavy-tailed Lévy processes [13, 31] to stochastic integrals driven by regularly varying Lévy processes. Note that a multivariate version of this result is also at hand.

**4. Proofs.** This section contains the proofs of the main results. For auxiliary results and technical lemmas, we refer to Section 5.

Throughout the rest of the paper we use the notation $B_{\mathbf{x},r}$ for the open ball in a metric space $(\mathbb{E}, \rho)$ with radius $r$, that is, $B_{\mathbf{x},r} = \{\mathbf{y} \in \mathbb{E} : \rho(\mathbf{y}, \mathbf{x}) < r\}$. The complement of a set $B \subset \mathbb{E}$ is denoted by $B^c$. The space $\mathbb{E}$ will usually be $\mathbb{D}$ or $\mathbb{R}^d$.

REMARK 4.1. The Portmanteau theorem implies that $\mathbf{X} \in \mathrm{RV}((a_n), m, \overline{\mathbb{D}}_0)$ if and only if $\limsup_{n \to \infty} n\,\mathrm{P}(\mathbf{X} \in a_n F) \leq m(F)$ and $\liminf_{n \to \infty} n\,\mathrm{P}(\mathbf{X} \in a_n G) \geq m(G)$ for all closed $F$ and open $G$ in $\mathbb{D}$ bounded away from $\mathbf{0}$. If there exist arbitrary small numbers $\delta > 0$ such that $\lim_{n \to \infty} n\,\mathrm{P}(|\mathbf{X}|_\infty \geq a_n \delta) = m(B_{\mathbf{0},\delta}^c)$, then it is straightforward to show that $\mathbf{X} \in \mathrm{RV}((a_n), m, \overline{\mathbb{D}}_0)$ if and only if $\limsup_{n \to \infty} n\,\mathrm{P}(\mathbf{X} \in a_n F) \leq m(F)$ for all closed $F$ in $\mathbb{D}$ bounded away from $\mathbf{0}$.

PROOF OF THEOREM 3.1. Take $B \in \mathcal{B}(\overline{\mathbb{D}}_0) \cap \mathbb{D}$, bounded away from $\mathbf{0}$, that is, $B \subset B_{\mathbf{0},\varepsilon}^c$ for some $\varepsilon > 0$, with $\mathrm{E}(m \circ \phi_{\mathbf{Y}}^{-1}(\partial B)) = 0$. By assumption, $\mathrm{E}(m(\mathrm{Disc}(\phi_{\mathbf{Y}}))) = 0$, and hence, there exists an $\Omega_0 \in \mathcal{F}$ with $\mathrm{P}(\Omega_0) = 1$ such that $m(\mathrm{Disc}(\phi_{\mathbf{Y}(\omega)})) = 0$ and $\mathbf{Y}(\omega) \neq \mathbf{0}$ for $\omega \in \Omega_0$. Let $d_B$ denote the shortest distance to the set $B$: $d_B = \inf\{|\mathbf{x}|_\infty : \mathbf{x} \in B\}$. We have

$$n\,\mathrm{P}(a_n^{-1} \mathbf{Y}\mathbf{X} \in B)$$

$$= n\,\mathrm{P}(a_n^{-1}\mathbf{Y}1_{(0,M)}(|\mathbf{Y}|_\infty)\mathbf{X} \in B) + n\,\mathrm{P}(a_n^{-1}\mathbf{Y}1_{[M,\infty)}(|\mathbf{Y}|_\infty)\mathbf{X} \in B)$$

$$= \int_{\{0<|\mathbf{y}|_\infty<M\}} \underbrace{n\,\mathrm{P}(a_n^{-1}\mathbf{y}\mathbf{X} \in B)}_{f_n(\mathbf{y})} \mathrm{P}(\mathbf{Y} \in d\mathbf{y})$$

$$+ n\,\mathrm{P}(a_n^{-1}\mathbf{Y}1_{[M,\infty)}(|\mathbf{Y}|_\infty)\mathbf{X} \in B).$$

Applying Theorem 2.1 yields $\lim_{n\to\infty} f_n(\mathbf{y}) = m \circ \phi_\mathbf{y}^{-1}(B)$ for each $\mathbf{y} \neq \mathbf{0}$. We want to show that

(4.1) $$\lim_{n\to\infty} \int_{\{0<|\mathbf{y}|_\infty<M\}} f_n(\mathbf{y})\,\mathrm{P}(\mathbf{Y} \in d\mathbf{y}) = \mathrm{E}(1_{(0,M)}(|\mathbf{Y}|_\infty) m \circ \phi_\mathbf{Y}^{-1}(B)),$$

(4.2) $$\limsup_{n\to\infty} n\,\mathrm{P}(a_n^{-1}\mathbf{Y}1_{[M,\infty)}(|\mathbf{Y}|_\infty)\mathbf{X} \in B) \leq C(M), \qquad \lim_{M\to\infty} C(M) = 0,$$

from which the conclusion follows by letting $M \to \infty$. To show (4.1), we use Pratt's theorem ([24], Theorem 1). For $0 < |\mathbf{y}|_\infty < M$,

$$f_n(\mathbf{y}) \leq n\,\mathrm{P}(a_n^{-1}|\mathbf{X}|_\infty > d_B/M) = G_n,$$

where $\lim_{n\to\infty} G_n = G = m(B_{\mathbf{0},d_B/M}^c) < \infty$. Clearly, as $n \to \infty$,

$$\int_{\{0<|\mathbf{y}|_\infty<M\}} G_n\,\mathrm{P}(\mathbf{Y} \in d\mathbf{y}) = \mathrm{P}(|\mathbf{Y}|_\infty < M)G_n \to \mathrm{P}(|\mathbf{Y}|_\infty < M)G.$$

Hence, Pratt's theorem can be applied from which follows that (4.1) holds. It remains to show (4.2). Applying Breiman's result (3.1) yields

$$\limsup_{n\to\infty} n\,\mathrm{P}(a_n^{-1}\mathbf{Y}1_{[M,\infty)}(|\mathbf{Y}|_\infty)\mathbf{X} \in B)$$

$$\leq \limsup_{n\to\infty} n\,\mathrm{P}(|\mathbf{Y}|_\infty 1_{[M,\infty)}(|\mathbf{Y}|_\infty)|\mathbf{X}|_\infty > a_n d_B)$$

$$= \mathrm{E}(|\mathbf{Y}|_\infty^\alpha 1_{[M,\infty)}(|\mathbf{Y}|_\infty)) m(B_{\mathbf{0},d_B}^c).$$

Since $\mathrm{E}(|\mathbf{Y}|_\infty^\alpha) < \infty$, it follows that $\lim_{M\to\infty} \mathrm{E}(|\mathbf{Y}|_\infty^\alpha 1_{(M,\infty)}(|\mathbf{Y}|_\infty)) = 0$. This proves (4.2). Thus, we have shown that

$$\limsup_{n\to\infty} n\,\mathrm{P}(a_n^{-1}\mathbf{Y}\mathbf{X} \in B) \leq \mathrm{E}(1_{[0,M)}(|\mathbf{Y}|_\infty) m \circ \phi_\mathbf{Y}^{-1}(B))$$

$$+ \mathrm{E}(|\mathbf{Y}|_\infty^\alpha 1_{[M,\infty)}(|\mathbf{Y}|_\infty)) m(B_{\mathbf{0},d_B}^c),$$

$$\liminf_{n\to\infty} n\,\mathrm{P}(a_n^{-1}\mathbf{Y}\mathbf{X} \in B) \geq \mathrm{E}(1_{[0,M)}(|\mathbf{Y}|_\infty) m \circ \phi_\mathbf{Y}^{-1}(B)).$$

Letting $M \to \infty$ now yields

$$\lim_{n\to\infty} n\,\mathrm{P}(a_n^{-1}\mathbf{Y}\mathbf{X} \in B) = \mathrm{E}(m \circ \phi_\mathbf{Y}^{-1}(B)).$$

Since $m(\overline{\mathbb{D}}_0 \setminus \mathbb{D}) = 0$ and $B \in \mathcal{B}(\overline{\mathbb{D}}_0) \cap \mathbb{D}$ with $\mathrm{E}(m \circ \phi_\mathbf{Y}^{-1}(\partial B)) = 0$ was arbitrary, the conclusion follows. □

PROOF OF THEOREM 3.2. Take $\varepsilon > 0$ and a closed set $F \in \mathcal{B}(\mathbb{D})$ with $d^\circ(\mathbf{0}, F) = \inf_{\mathbf{z} \in F} d^\circ(\mathbf{0}, \mathbf{z}) > \varepsilon$. Define $F_\varepsilon = \{\mathbf{x} \in \mathbb{D} : d^\circ(\mathbf{x}, F) \leq \varepsilon\}$. Then $F, F_\varepsilon \in \mathcal{B}(\overline{\mathbb{D}}_0)$ and both $F$ and $F_\varepsilon$ are closed and bounded in $\overline{\mathbb{D}}_0$. Take $\delta > \varepsilon$. Notice that

$$P(d^\circ(a_n^{-1}\mathbf{X}, a_n^{-1}\mathbf{Y}) \leq \varepsilon \mid |\mathbf{Y}|_\infty > a_n\delta) = \frac{P(d^\circ(a_n^{-1}\mathbf{X}, a_n^{-1}\mathbf{Y}) \leq \varepsilon, |\mathbf{Y}|_\infty > a_n\delta)}{P(|\mathbf{Y}|_\infty > a_n\delta)}$$

$$\leq \frac{P(|\mathbf{X}|_\infty > a_n(\delta - \varepsilon))}{P(|\mathbf{Y}|_\infty > a_n\delta)}.$$

Hence, the first part of (3.3) yields

$$\limsup_{n \to \infty} n P(|\mathbf{Y}|_\infty > a_n\delta) \leq \limsup_{n \to \infty} \frac{n P(|\mathbf{X}|_\infty > a_n(\delta - \varepsilon))}{P(d^\circ(a_n^{-1}\mathbf{X}, a_n^{-1}\mathbf{Y}) \leq \varepsilon \mid |\mathbf{Y}|_\infty > a_n\delta)}$$

$$= (\delta - \varepsilon)^{-\alpha} m(B_{\mathbf{0},1}^c) \in (0, \infty).$$

Similarly, switching from $\mathbf{Y}$ to $\mathbf{X}$ in the second to last expression above and applying the second part of (3.3) we obtain

$$\liminf_{n \to \infty} n P(|\mathbf{Y}|_\infty > a_n\delta) \geq (\delta + \varepsilon)^{-\alpha} m(B_{\mathbf{0},1}^c).$$

Since $\varepsilon$ may be chosen arbitrarily small we conclude that

$$\lim_{n \to \infty} n P(|\mathbf{Y}|_\infty \geq a_n\delta) = \lim_{n \to \infty} n P(|\mathbf{Y}|_\infty > a_n\delta) = \delta^{-\alpha} m(B_{\mathbf{0},1}^c).$$

Moreover, we observe that for $\delta \in (\varepsilon, d^\circ(\mathbf{0}, F))$

$$n P(a_n^{-1}\mathbf{Y} \in F) \leq n P(a_n^{-1}\mathbf{Y} \in F, d^\circ(a_n^{-1}\mathbf{X}, a_n^{-1}\mathbf{Y}) \geq \varepsilon) + n P(a_n^{-1}\mathbf{X} \in F_\varepsilon)$$

$$\leq n P(|\mathbf{Y}|_\infty > a_n\delta) P(d^\circ(a_n^{-1}\mathbf{X}, a_n^{-1}\mathbf{Y}) \geq \varepsilon \mid |\mathbf{Y}|_\infty > a_n\delta)$$

$$+ n P(a_n^{-1}\mathbf{X} \in F_\varepsilon).$$

Since $F_\varepsilon$ is closed the hypotheses and the Portmanteau theorem imply that

$$\limsup_{n \to \infty} n P(a_n^{-1}\mathbf{Y} \in F) \leq \limsup_{n \to \infty} n P(a_n^{-1}\mathbf{X} \in F_\varepsilon) \leq m(F_\varepsilon).$$

Since $F$ is closed, $F_\varepsilon \downarrow F$ as $\varepsilon \downarrow 0$. Hence, $\limsup_{n \to \infty} n P(a_n^{-1}\mathbf{Y} \in F) \leq m(F)$ and the conclusion follows from the Portmanteau theorem and Remark 4.1. □

PROOF OF THEOREM 3.3. For $\gamma > 0$, we say that an element $\mathbf{x} \in \mathbb{D}$ has $\gamma$-oscillation $p$ times in $[0, 1]$ if there exist $0 \leq t_0 < t_1 < \cdots < t_p \leq 1$ such that $|\mathbf{x}_{t_i} - \mathbf{x}_{t_{i-1}}| > \gamma$ for each $i = 1, \ldots, p$. We write

$$B(p, \gamma, [0, 1]) = \{\mathbf{x} \in \mathbb{D} : \mathbf{x} \text{ has } \gamma\text{-oscillation } p \text{ times in } [0, 1]\}.$$

We start with the first claim. Take $\varepsilon > 0$ and set $\tau = \tau(\mathbf{X})$. Since

$$\mathrm{P}(d^\circ(a_n^{-1}\mathbf{X}, a_n^{-1}\Delta\mathbf{X}_\tau 1_{[\tau,1]}) > \varepsilon \mid |\mathbf{X}|_\infty > a_n)$$
$$= \frac{n\,\mathrm{P}(d^\circ(a_n^{-1}\mathbf{X}, a_n^{-1}\Delta\mathbf{X}_\tau 1_{[\tau,1]}) > \varepsilon, |\mathbf{X}|_\infty > a_n)}{n\,\mathrm{P}(|\mathbf{X}|_\infty > a_n)},$$

it is sufficient to show that the numerator tends to zero as $n \to \infty$. Moreover, we can, without loss of generality, take $\varepsilon \leq 1$. We have that

$$n\,\mathrm{P}(d^\circ(a_n^{-1}\mathbf{X}, a_n^{-1}\Delta\mathbf{X}_\tau 1_{[\tau,1]}) > \varepsilon, |a_n^{-1}\mathbf{X}|_\infty > 1)$$
$$= n\,\mathrm{P}(d^\circ(a_n^{-1}\mathbf{X}, a_n^{-1}\Delta\mathbf{X}_\tau 1_{[\tau,1]}) > \varepsilon, |\mathbf{X}|_\infty > a_n, a_n^{-1}\mathbf{X} \notin B(2, \varepsilon/4, [0,1]))$$
$$+ n\,\mathrm{P}(d^\circ(a_n^{-1}\mathbf{X}, a_n^{-1}\Delta\mathbf{X}_\tau 1_{[\tau,1]}) > \varepsilon, |\mathbf{X}|_\infty > a_n, a_n^{-1}\mathbf{X} \in B(2, \varepsilon/4, [0,1]))$$
$$= p_n + q_n.$$

Note that

$$q_n \leq n\,\mathrm{P}(a_n^{-1}\mathbf{X} \in B(2, \varepsilon/4, [0,1])) \to 0,$$

by Lemma 21 in [17]. Note also that if $\mathbf{x} \in \mathbb{D}$, $\mathbf{x}_0 = \mathbf{0}$ and $\mathbf{x} \in B_{\mathbf{0},1}^c \cap B(2, \varepsilon/4, [0,1])^c$, then there exists $t_0 \in (0,1)$ such that

$$\mathbf{x} \notin B(1, \varepsilon/4, [0, t_0)), \qquad |\Delta\mathbf{x}_{t_0}| > \varepsilon/2 \quad \text{and} \quad \mathbf{x} \notin B(1, \varepsilon/4, [t_0, 1]).$$

It follows that $t_0 = \tau(\mathbf{x})$, $d^\circ(\mathbf{x}, \Delta\mathbf{x}_{t_0} 1_{[t_0,1]}) < \varepsilon$, and hence, that

$$p_n \leq n\,\mathrm{P}(d^\circ(a_n^{-1}\mathbf{X}, a_n^{-1}\Delta\mathbf{X}_\tau 1_{[\tau,1]}) > \varepsilon, d^\circ(a_n^{-1}\mathbf{X}, a_n^{-1}\Delta\mathbf{X}_\tau 1_{[\tau,1]}) < \varepsilon) = 0.$$

This completes the proof of the first claim. For the second claim, take w.l.g. $\varepsilon \in (0,1)$ and note that, as $n \to \infty$,

$$n\,\mathrm{P}(d^\circ(a_n^{-1}\mathbf{X}, a_n^{-1}\Delta\mathbf{X}_\tau 1_{[\tau,1]}) > \varepsilon, |\Delta\mathbf{X}_\tau 1_{[\tau,1]}|_\infty > a_n)$$
$$\leq n\,\mathrm{P}(a_n^{-1}\mathbf{X} \in B(2, \varepsilon/2, [0,1])) \to 0.$$

By Lemma 2.1, $\mathbf{X} \in \mathrm{RV}((a_n), m, \overline{\mathbb{D}}_0)$ implies $\mathbf{X}_1 \in \mathrm{RV}((a_n), m_1, \overline{\mathbb{R}}_0^d)$ and

$$\lim_{n \to \infty} n\,\mathrm{P}(a_n^{-1}\mathbf{X}_1 \in B_{\mathbf{0},1}^c) = \lim_{n \to \infty} n\nu(a_n B_{\mathbf{0},1}^c) = m(B_{\mathbf{0},1}^c) > 0,$$

where $\nu$ denotes the Lévy measure of $\mathbf{X}_1$. Hence, as $n \to \infty$,

$$n\,\mathrm{P}(|\Delta\mathbf{X}_\tau 1_{[\tau,1]}|_\infty > a_n) = n\,\mathrm{P}(\xi([0,1] \times a_n B_{\mathbf{0},1}^c) > 0)$$
$$= n(1 - e^{-\nu(a_n B_{\mathbf{0},1}^c)})$$
$$\sim n\nu(a_n B_{\mathbf{0},1}^c) \to m(B_{\mathbf{0},1}^c) > 0.$$

It follows that $\lim_{n \to \infty} \mathrm{P}(d^\circ(a_n^{-1}\mathbf{X}, a_n^{-1}\Delta\mathbf{X}_\tau 1_{[\tau,1]}) > \varepsilon \mid |\Delta\mathbf{X}_\tau 1_{[\tau,1]}|_\infty > a_n) = 0$, and the conclusion follows from Theorem 3.2. □

PROOF OF THEOREM 3.4. As usual we set $\tau = \tau(\mathbf{X})$. The outline of the proof is as follows:

(i) Show that
$$\lim_{n\to\infty} n\,P(d^\circ(a_n^{-1}(\mathbf{Y}\cdot\mathbf{X}), a_n^{-1}\mathbf{Y}_\tau\Delta\mathbf{X}_\tau 1_{[\tau,1]}) > \varepsilon, |(\mathbf{Y}\cdot\mathbf{X})|_\infty > a_n) = 0.$$

(ii) Show that, for each $\delta > 0$, $\lim_{n\to\infty} n\,P(|(\mathbf{Y}\cdot\mathbf{X})|_\infty > a_n\delta) = m^*(B_{0,\delta}^c)$.
From (i) and (ii) we conclude that (3.6) holds.

(iii) Show that $\mathbf{Y}_\tau\Delta\mathbf{X}_\tau 1_{[\tau,1]} \in \mathrm{RV}((a_n), m^*, \overline{\mathbb{D}}_0)$.

Then (3.7) holds because for each $\varepsilon > 0$

$$1 \geq \lim_{n\to\infty} P(d^\circ(a_n^{-1}(\mathbf{Y}\cdot\mathbf{X}), a_n^{-1}\mathbf{Y}_\tau\Delta\mathbf{X}_\tau 1_{[\tau,1]}) \leq \varepsilon \mid |\mathbf{Y}_\tau\Delta\mathbf{X}_\tau| > a_n)$$

$$= \lim_{n\to\infty} \frac{P(d^\circ(a_n^{-1}(\mathbf{Y}\cdot\mathbf{X}), a_n^{-1}\mathbf{Y}_\tau\Delta\mathbf{X}_\tau 1_{[\tau,1]}) \leq \varepsilon, |\mathbf{Y}_\tau\Delta\mathbf{X}_\tau| > a_n)}{P(|\mathbf{Y}_\tau\Delta\mathbf{X}_\tau| > a_n)}$$

$$\geq \lim_{n\to\infty} \frac{P(d^\circ(a_n^{-1}(\mathbf{Y}\cdot\mathbf{X}), a_n^{-1}\mathbf{Y}_\tau\Delta\mathbf{X}_\tau 1_{[\tau,1]}) \leq \varepsilon, |(\mathbf{Y}\cdot\mathbf{X})|_\infty > a_n(1+\varepsilon))}{P(|(\mathbf{Y}\cdot\mathbf{X})|_\infty > a_n(1+\varepsilon))}$$

$$\times \frac{P(|(\mathbf{Y}\cdot\mathbf{X})|_\infty > a_n(1+\varepsilon))}{P(|\mathbf{Y}_\tau\Delta\mathbf{X}_\tau| > a_n)}$$

$$= (1+\varepsilon)^{-\alpha}.$$

Finally, Theorem 3.2 gives the conclusion.

(i) Take $\varepsilon > 0$; w.l.g., we can take $\varepsilon \leq 1$. Then, writing $\mathbf{X} = \widetilde{\mathbf{X}} + \mathbf{J}$, we have

$$\{d^\circ(a_n^{-1}(\mathbf{Y}\cdot\mathbf{X}), a_n^{-1}\mathbf{Y}_\tau\Delta\mathbf{X}_\tau 1_{[\tau,1]}) > \varepsilon, |(\mathbf{Y}\cdot\mathbf{X})|_\infty > a_n\}$$

$$\subset \{d^\circ(a_n^{-1}(\mathbf{Y}\cdot\mathbf{X}), a_n^{-1}\mathbf{Y}_\tau\Delta\mathbf{X}_\tau 1_{[\tau,1]}) > \varepsilon, |(\mathbf{Y}\cdot\mathbf{X})|_\infty > a_n,$$
$$|(\mathbf{Y}\cdot\widetilde{\mathbf{X}})|_\infty > a_n\varepsilon/2\}$$

$$\cup \{d^\circ(a_n^{-1}(\mathbf{Y}\cdot\mathbf{X}), a_n^{-1}\mathbf{Y}_\tau\Delta\mathbf{X}_\tau 1_{[\tau,1]}) > \varepsilon, |(\mathbf{Y}\cdot\mathbf{X})|_\infty > a_n,$$
$$|(\mathbf{Y}\cdot\widetilde{\mathbf{X}})|_\infty \leq a_n\varepsilon/2\}$$

$$= A_n \cup B_n.$$

We will show that $\lim_{n\to\infty} n\,P(A_n) = 0$ and $\lim_{n\to\infty} n\,P(B_n) = 0$. Note that $A_n \subset \{|(\mathbf{Y}\cdot\widetilde{\mathbf{X}})|_\infty > a_n\varepsilon/2\}$. By a standard regular variation argument, $\mathbf{X}_1 \in \mathrm{RV}_\alpha((a_n), \mu, \overline{\mathbb{R}}_0^d)$ implies that the sequence $(a_n)$ is regularly varying with index $1/\alpha$. By construction, $\widetilde{\mathbf{X}}$ is a Lévy process with bounded jumps, so Lemma 5.5 gives

(4.3) $$\lim_{n\to\infty} n\,P(A_n) \leq \lim_{n\to\infty} n\,P(|(\mathbf{Y}\cdot\widetilde{\mathbf{X}})|_\infty > a_n\varepsilon/2) = 0.$$

Next we consider $\lim_{n\to\infty} n\,P(B_n)$. First we note that for any $\varepsilon \in (0,1)$ and $\mathbf{x}, \mathbf{y} \in \mathbb{D}$, $|\mathbf{x}+\mathbf{y}|_\infty > 1$ and $|\mathbf{y}|_\infty \leq \varepsilon/2$ implies $|\mathbf{x}|_\infty > 1/2$. Hence,

$$B_n \subset \{d^\circ(a_n^{-1}(\mathbf{Y}\cdot\mathbf{J}), a_n^{-1}\mathbf{Y}_\tau\Delta\mathbf{X}_\tau 1_{[\tau,1]}) > \varepsilon/2, |(\mathbf{Y}\cdot\mathbf{X})|_\infty > a_n,$$

$$|(\mathbf{Y}\cdot\tilde{\mathbf{X}})|_\infty \leq a_n\varepsilon/2\}$$
$$\subset \{d^\circ(a_n^{-1}(\mathbf{Y}\cdot\mathbf{J}), a_n^{-1}\mathbf{Y}_\tau\Delta\mathbf{X}_\tau 1_{[\tau,1]}) > \varepsilon/2, |(\mathbf{Y}\cdot\mathbf{J})|_\infty > a_n/2\} = C_n.$$

Let $((\tau_k, \mathbf{Z}_k))_{k\geq 1}$ be the points of the PRM $\xi(\cdot \cap \{[0,1] \times B_{\mathbf{0},1}^c\})$ [see (2.3)] and note that $(\mathbf{Z}_k)$ is an i.i.d. sequence with $\mathbf{Z}_1 \in \mathrm{RV}((a_n), m_1, \overline{\mathbb{R}}_0^d)$. We have

$$(\mathbf{Y}\cdot\mathbf{J})_t = \sum_{k=1}^{N_t} \mathbf{Y}_{\tau_k}\mathbf{Z}_k,$$

where $N_t = \xi((0,t] \times B_{\mathbf{0},1}^c)$. Note that $(N_t)$ and $(\mathbf{Z}_k)$ are independent and, since $\mathbf{Y}$ is predictable, for every $k$, $\mathbf{Y}_{\tau_k}$ and $\mathbf{Z}_k$ are independent. For $\beta \in (1/2, 1)$, let

$$\mathbf{J}_n = \sum_{k=1}^{N_1} \mathbf{Z}_k 1_{(a_n^\beta, \infty)}(|\mathbf{Z}_k|)1_{[\tau_k,1]},$$

that is, $\mathbf{J}_n$ consists of the jumps with norm larger than $a_n^\beta$. Then

$$C_n \subset \{d^\circ(a_n^{-1}(\mathbf{Y}\cdot\mathbf{J}), a_n^{-1}\mathbf{Y}_\tau\Delta\mathbf{X}_\tau 1_{[\tau,1]}) > \varepsilon/2, |(\mathbf{Y}\cdot(\mathbf{J}-\mathbf{J}_n))|_\infty > a_n/4\}$$
$$\cup \{d^\circ(a_n^{-1}(\mathbf{Y}\cdot\mathbf{J}), a_n^{-1}\mathbf{Y}_\tau\Delta\mathbf{X}_\tau 1_{[\tau,1]}) > \varepsilon/2, |(\mathbf{Y}\cdot\mathbf{J}_n)|_\infty > a_n/4\}$$
$$\subset \{|(\mathbf{Y}\cdot(\mathbf{J}-\mathbf{J}_n))|_\infty > a_n/4\}$$
$$\cup \underbrace{\{d^\circ(a_n^{-1}(\mathbf{Y}\cdot\mathbf{J}), a_n^{-1}\mathbf{Y}_\tau\Delta\mathbf{X}_\tau 1_{[\tau,1]}) > \varepsilon/2, |(\mathbf{Y}\cdot\mathbf{J}_n)|_\infty > a_n/4\}}_{D_n}.$$

Introduce $M_n = \sum_{k=1}^{N_1} 1_{(a_n^\beta, \infty)}(|\mathbf{Z}_k|)$, the number of jumps with norm larger than $a_n^\beta$, and note that on $\{M_n = 1\}$, we have $\Delta\mathbf{X}_\tau 1_{[\tau,1]} = \mathbf{J}_n$. Hence,

$$D_n \subset \{d^\circ(a_n^{-1}(\mathbf{Y}\cdot\mathbf{J}), a_n^{-1}(\mathbf{Y}\cdot\mathbf{J}_n)) > \varepsilon/2, M_n = 1\} \cup \{M_n \geq 2\}$$
$$\subset \{|(\mathbf{Y}\cdot(\mathbf{J}-\mathbf{J}_n))|_\infty > a_n\varepsilon/2\} \cup \{M_n \geq 2\}.$$

Putting everything together, we see that, with $\delta < \min(\varepsilon/2, 1/4)$, the set $B_n$ satisfies

$$n\mathrm{P}(B_n) \leq 2n\mathrm{P}(|(\mathbf{Y}\cdot(\mathbf{J}-\mathbf{J}_n))|_\infty > a_n\delta) + n\mathrm{P}(M_n \geq 2)$$
$$\leq 2n\mathrm{P}\left(\sum_{k=0}^{N_1} |\mathbf{Y}_{\tau_k}||\mathbf{Z}_k|1_{[0,a_n^\beta]}(|\mathbf{Z}_k|) > a_n\delta\right) + n\mathrm{P}(M_n \geq 2).$$

The first term converges to zero by Lemma 5.3 and Remark 5.1 and for the second term, we apply Lemma 5.4. This proves $\lim_{n\to\infty} n\mathrm{P}(B_n) = 0$ and, hence, we have shown (i).

(ii) Take $\delta > 0$ and note that $m^*(\partial B_{0,\delta}^c) = 0$. Using (4.3) and applying Proposition 5.1 we find that for each $\varepsilon \in (0, \delta)$

$$\liminf_{n \to \infty} n \, \mathrm{P}(|(\mathbf{Y} \cdot \mathbf{X})|_\infty > a_n \delta)$$
$$= \liminf_{n \to \infty} n \, \mathrm{P}(|(\mathbf{Y} \cdot \mathbf{X})|_\infty > a_n \delta, |(\mathbf{Y} \cdot \widetilde{\mathbf{X}})|_\infty > a_n \varepsilon)$$
$$+ \liminf_{n \to \infty} n \, \mathrm{P}(|(\mathbf{Y} \cdot \mathbf{X})|_\infty > a_n \delta, |(\mathbf{Y} \cdot \widetilde{\mathbf{X}})|_\infty \leq a_n \varepsilon)$$
$$\geq \liminf_{n \to \infty} n \, \mathrm{P}(|(\mathbf{Y} \cdot \mathbf{J})|_\infty > a_n(\delta + \varepsilon), |(\mathbf{Y} \cdot \widetilde{\mathbf{X}})|_\infty \leq a_n \varepsilon)$$
$$= \liminf_{n \to \infty} n \, \mathrm{P}(|(\mathbf{Y} \cdot \mathbf{J})|_\infty > a_n(\delta + \varepsilon))$$
$$= m^*(B_{\mathbf{0}, \delta+\varepsilon}^c).$$

Similarly,

$$\limsup_{n \to \infty} n \, \mathrm{P}(|(\mathbf{Y} \cdot \mathbf{X})|_\infty > a_n \delta)$$
$$= \limsup_{n \to \infty} n \, \mathrm{P}(|(\mathbf{Y} \cdot \mathbf{X})|_\infty > a_n \delta, |(\mathbf{Y} \cdot \widetilde{\mathbf{X}})|_\infty > a_n \varepsilon)$$
$$+ \limsup_{n \to \infty} n \, \mathrm{P}(|(\mathbf{Y} \cdot \mathbf{X})|_\infty > a_n \delta, |(\mathbf{Y} \cdot \widetilde{\mathbf{X}})|_\infty \leq a_n \varepsilon)$$
$$\leq \limsup_{n \to \infty} n \, \mathrm{P}(|(\mathbf{Y} \cdot \mathbf{J})|_\infty > a_n(\delta - \varepsilon))$$
$$= m^*(B_{\mathbf{0}, \delta-\varepsilon}^c).$$

Then (ii) follows by letting $\varepsilon \to 0$.

(iii) Take closed $B \in \mathcal{B}(\overline{\mathbb{D}}_0)$ bounded away from $\mathbf{0}$ and set $d_B = \inf\{|\mathbf{x}|_\infty : \mathbf{x} \in B\}$. We will show that

(4.4) $$\limsup_{n \to \infty} n \, \mathrm{P}(a_n^{-1} \mathbf{Y}_\tau \Delta \mathbf{X}_\tau 1_{[\tau,1]} \in B) \leq m^*(B)$$

and for each $\delta > 0$

(4.5) $$\lim_{n \to \infty} n \, \mathrm{P}(|\mathbf{Y}_\tau \Delta \mathbf{X}_\tau| > a_n \delta) = m^*(B_{0,\delta}^c).$$

Then, (iii) follows from the Portmanteau theorem.

For an element $\mathbf{z} \in \mathbb{D}$, we denote by $S(\mathbf{z}) = \Delta \mathbf{z}_{\tau(\mathbf{z})} 1_{[\tau(\mathbf{z}),1]}$ the step function with one step at $\tau(\mathbf{z})$. Note that $S(\mathbf{X})$ equals either $S(\mathbf{J})$ (if $|\Delta \mathbf{X}_\tau| \geq 1$) or $S(\widetilde{\mathbf{X}})$, so

$$\{\mathbf{Y}_\tau \Delta \mathbf{X}_\tau 1_{[\tau,1]} \in a_n B\} = \{(\mathbf{Y} \cdot S(\mathbf{J})) \in a_n B\} \cup \{(\mathbf{Y} \cdot S(\widetilde{\mathbf{X}})) \in a_n B\}.$$

Since the jumps of $\widetilde{\mathbf{X}}$ are bounded by 1 and the sequence $(a_n)$ is regularly varying with index $1/\alpha$, we have

$$n \, \mathrm{P}((\mathbf{Y} \cdot S(\widetilde{\mathbf{X}})) \in a_n B) \leq n \, \mathrm{P}(|\mathbf{Y}|_\infty > a_n d_B) \leq n(a_n d_B)^{-\alpha-\delta} \, \mathrm{E}(|\mathbf{Y}|_\infty^{\alpha+\delta}) \to 0.$$

For the term involving $S(\mathbf{J})$, we write

$$\{(\mathbf{Y} \cdot S(\mathbf{J})) \in a_n B\} = \{(\mathbf{Y} \cdot S(\mathbf{J}_n)) \in a_n B\} \cup \{(\mathbf{Y} \cdot S(\mathbf{J} - \mathbf{J}_n)) \in a_n B\}.$$

By Lemma 5.3 and Remark 5.1,

$$n \operatorname{P}((\mathbf{Y} \cdot S(\mathbf{J} - \mathbf{J}_n)) \in a_n B) \leq n \operatorname{P}\left(\sum_{k=1}^{N_1} |\mathbf{Y}_{\tau_k}||\mathbf{Z}_k| 1_{[1,a_n^\beta]}(|\mathbf{Z}_k|) > a_n d_B\right) \to 0.$$

Moreover,

$$(4.6) \quad n \operatorname{P}((\mathbf{Y} \cdot S(\mathbf{J}_n)) \in a_n B) = n \operatorname{P}((\mathbf{Y} \cdot S(\mathbf{J}_n)) \in a_n B, M_n = 1)$$
$$(4.7) \qquad\qquad\qquad\qquad + n \operatorname{P}((\mathbf{Y} \cdot S(\mathbf{J}_n)) \in a_n B, M_n \geq 2).$$

The term in (4.7) is less than or equal to $n \operatorname{P}(M_n \geq 2)$, which converges to 0 by Lemma 5.4. On $\{M_n = 1\}$ we have $S(\mathbf{J}_n) = \mathbf{J}_n$ so the term (4.6) satisfies

$$\limsup_{n \to \infty} n \operatorname{P}((\mathbf{Y} \cdot S(\mathbf{J}_n)) \in a_n B, M_n = 1)$$
$$= \limsup_{n \to \infty} n \operatorname{P}((\mathbf{Y} \cdot \mathbf{J}_n) \in a_n B, M_n = 1)$$
$$\leq m^*(B),$$

by following the lines of the proof of Proposition 5.1. This proves (4.4). The proof of (4.5) is similar;

$$n \operatorname{P}(|\mathbf{Y}_\tau \Delta \mathbf{X}_\tau| > a_n \delta) \sim n \operatorname{P}(|(\mathbf{Y} \cdot \mathbf{J}_n)|_\infty > a_n \delta, M_n = 1) \to m^*(B_{0,\delta}^c)$$

as $n \to \infty$. This completes the proof. $\square$

## 5. Auxiliary results.

LEMMA 5.1. *Let $\mathbf{X}$ and $\mathbf{Y}$ be stochastic processes with sample paths in $\mathbb{D}$, with $\mathbf{X}$ being a Lévy process satisfying $\mathbf{X} \in \operatorname{RV}((a_n), m, \overline{\mathbb{D}}_0)$. Then $\operatorname{E}(m(\operatorname{Disc}(\phi_\mathbf{Y}))) = 0$.*

PROOF. For $B \in \mathcal{B}(\overline{\mathbb{R}}_0^d) \cap \mathbb{R}^d$ and $T \in \mathcal{B}([0,1))$, denote

$$A_{B,T} = \{\mathbf{x} = \mathbf{y} 1_{[v,1]} : \mathbf{y} \in B, v \in T\}.$$

By the representation (2.5) of $m$, we have $m(A_{B,T}) = \mu(B)\lambda(T)$, where $\mu$ is the limit measure of $\mathbf{X}_1$ and $\lambda$ denotes the Lebesgue measure on $[0,1)$. Recall the notation $\mathcal{V} \subset \mathbb{D}$ from (2.4) for the support of $m$, which is the set of step functions with one step. Take an arbitrary $\omega \in \Omega$ and let $D_{\mathbf{Y}(\omega)}$ be the discontinuity points of $\mathbf{Y}(\omega)$. Then

$$\operatorname{Disc}(\phi_{\mathbf{Y}(\omega)}) \cap \mathcal{V} = \bigcup_{\varepsilon \in (0,\infty) \cap \mathbb{Q}} A_{B_{\mathbf{0},\varepsilon}^c, D_{\mathbf{Y}(\omega)}}.$$

Since $\mathbf{Y}(\omega) \in \mathbb{D}$, it follows that $D_{\mathbf{Y}(\omega)}$ is at most countable and $\lambda(D_{\mathbf{Y}(\omega)}) = 0$. Hence,

$$m(\mathrm{Disc}(\phi_{\mathbf{Y}(\omega)})) = m(\mathrm{Disc}(\phi_{\mathbf{Y}(\omega)}) \cap \mathcal{V})$$
$$\leq \sum_{\varepsilon \in (0,\infty) \cap \mathbb{Q}} m(A_{B_{\mathbf{0},\varepsilon}^c, D_{\mathbf{Y}(\omega)}}) = 0.$$

Since $\omega \in \Omega$ was arbitrary, we see that $\mathrm{E}(m(\mathrm{Disc}(\phi_\mathbf{Y}))) = 0$. □

LEMMA 5.2. *Let $(Z_k)$ be an i.i.d. sequence of nonnegative random variables, let $N$ be an $\mathbb{N}$-valued random variable and let $(Y_k)$ be a sequence of non-negative random variables. Suppose further that $(\mathcal{F}_k)$ is a filtration such that $Y_k$ is $\mathcal{F}_k$-measurable, $Z_k$ is $\mathcal{F}_{k+1}$-measurable and independent of $\mathcal{F}_k$ and $N$. Then, for each $x > 0$,*

$$(5.1) \qquad \mathrm{P}\left(\sum_{k=1}^N Y_k Z_k > x\right) \leq 2 \mathrm{P}\left(N \bigvee_{k=1}^N Y_k \widetilde{Z}_k > x\right),$$

*where $(\widetilde{Z}_k) \stackrel{d}{=} (Z_k)$ (possibly on an extended probability space) and $(\widetilde{Z}_k)$ is independent of $(Y_k)$ and $N$.*

PROOF. Let $\mathcal{F}_k' = \sigma(\mathcal{F}_k, \widetilde{Z}_1, \ldots, \widetilde{Z}_{k-1})$. By the assumptions, we have

$$(5.2) \qquad \mathrm{P}(Y_k Z_k \in \cdot | \mathcal{F}_k', N = m) = \mathrm{P}(Y_k \widetilde{Z}_k \in \cdot | \mathcal{F}_k', N = m).$$

Conditioning on $N$, we write

$$\mathrm{P}\left(\sum_{k=1}^N Y_k Z_k > x\right) = \sum_{m=1}^\infty \mathrm{P}\left(\sum_{k=1}^m Y_k Z_k > x \Big| N = m\right) \mathrm{P}(N = m)$$
$$\leq \sum_{m=1}^\infty \mathrm{P}\left(\bigvee_{k=1}^m Y_k Z_k > \frac{x}{m} \Big| N = m\right) \mathrm{P}(N = m).$$

Let $\tau = \min\{k : Y_k \widetilde{Z}_k > \frac{x}{m}\}$ and note that $\{\tau \geq k\} = \{\tau \leq k-1\}^c$ is $\mathcal{F}_k'$-measurable. Moreover,

$$\mathrm{P}\left(\bigvee_{k=1}^m Y_k Z_k > \frac{x}{m} \Big| N = m\right)$$
$$\leq \mathrm{P}\left(\bigvee_{k=1}^m Y_k Z_k > \frac{x}{m}, \tau \leq m \Big| N = m\right) + \mathrm{P}\left(\bigvee_{k=1}^m Y_k Z_k > \frac{x}{m}, \tau \geq m \Big| N = m\right)$$
$$\leq \mathrm{P}(\tau \leq m | N = m) + \sum_{k=1}^m \mathrm{P}\left(\tau \geq k, Y_k Z_k > \frac{x}{m} \Big| N = m\right).$$

Using (5.2), the last expression equals

$$P(\tau \leq m | N = m) + \sum_{k=1}^{m} P\left(\tau \geq k, Y_k \widetilde{Z}_k > \frac{x}{m} \Big| N = m\right)$$

$$= P(\tau \leq m | N = m) + \sum_{k=1}^{m} P(\tau = k | N = m)$$

$$= 2 P(\tau \leq m | N = m)$$

$$= 2 P\left(\bigvee_{k=1}^{m} Y_k \widetilde{Z}_k > \frac{x}{m} \Big| N = m\right).$$

Summing up over $m$, we arrive at (5.1), which proves the lemma. □

LEMMA 5.3. *Assume the hypotheses of Lemma 5.2. Suppose further that $Z_1 \in \mathrm{RV}_\alpha((a_n), \mu, (0, \infty])$ for some $\alpha > 0$ and that $\mathrm{E}(N^{\alpha+\gamma} \sum_{k=1}^{N} Y_k^{\alpha+\gamma}) < \infty$ for some $\gamma > 0$. Then, for every $\beta \in (0, 1)$,*

$$\lim_{n \to \infty} n P\left(\sum_{k=1}^{N} Y_k Z_k 1_{[0, a_n^\beta]}(Z_k) > a_n x\right) = 0, \qquad x > 0.$$

REMARK 5.1. Assume the hypotheses of Theorem 3.4, let $N$ be the number of jumps of **X** of norm greater than one and let $\tau_1, \ldots, \tau_N$ be the times of these jumps. Moreover, let $Y_k = |\mathbf{Y}_{\tau_k}|$, $\gamma = \delta/2$, $p = (\alpha + \delta)/(\alpha + \delta/2)$ and $q = (1 - 1/p)^{-1}$. Then, Lemma 5.3 applies. Indeed, using Hölder's inequality, we find that

$$E\left(N^{\alpha+\gamma} \sum_{k=1}^{N} Y_k^{\alpha+\gamma}\right) \leq E\left(N^{\alpha+\gamma+1} \bigvee_{k=1}^{N} Y_k^{\alpha+\gamma}\right)$$

$$\leq E(N^{q(\alpha+\gamma+1)})^{1/q} E(|\mathbf{Y}|_\infty^{p(\alpha+\gamma)})^{1/p}$$

$$= E(N^{q(\alpha+\gamma+1)})^{1/q} E(|\mathbf{Y}|_\infty^{\alpha+\delta})^{1/p} < \infty.$$

PROOF OF LEMMA 5.3. By Lemma 5.2,

$$n P\left(\sum_{k=1}^{N} Y_k Z_k 1_{[0, a_n^\beta]}(Z_k) > a_n x\right)$$

$$\leq 2n P\left(N \bigvee_{k=1}^{N} Y_k \widetilde{Z}_k 1_{[0, a_n^\beta]}(\widetilde{Z}_k) > a_n x\right).$$

Conditioning on $N$, we get

$$2n P\left(N \bigvee_{k=1}^{N} Y_k \widetilde{Z}_k 1_{[0, a_n^\beta]}(\widetilde{Z}_k) > a_n x\right)$$

$$= 2n \sum_{m=1}^{\infty} \mathrm{P}\left(m \bigvee_{k=1}^{m} Y_k \widetilde{Z}_k 1_{[0,a_n^\beta]}(\widetilde{Z}_k) > a_n x \Big| N = m\right) \mathrm{P}(N=m)$$

$$= 2n \sum_{m=1}^{\infty} \mathrm{P}\left(\bigcup_{k=1}^{m} \left\{Y_k \widetilde{Z}_k 1_{[0,a_n^\beta]}(\widetilde{Z}_k) > \frac{a_n x}{m} \Big| N = m\right\}\right) \mathrm{P}(N=m)$$

$$\leq 2n \sum_{m=1}^{\infty} \sum_{k=1}^{m} \mathrm{P}\left(Y_k \widetilde{Z}_k 1_{[0,a_n^\beta]}(\widetilde{Z}_k) > \frac{a_n x}{m} \Big| N = m\right) \mathrm{P}(N=m).$$

Denote the distribution of $\widetilde{Z}_k$ by $F$. By conditioning on $\widetilde{Z}_k$ and then using Markov's inequality, the last expression equals

$$2n \sum_{m=1}^{\infty} \sum_{k=1}^{m} \int_0^{a_n^\beta} \mathrm{P}\left(Y_k > \frac{a_n x}{mz} \Big| N = m\right) F(dz) \mathrm{P}(N=m)$$

$$\leq 2n \sum_{m=1}^{\infty} \sum_{k=1}^{m} \int_0^{a_n^\beta} \left(\frac{a_n x}{mz}\right)^{-(\alpha+\gamma)} E(Y_k^{\alpha+\gamma} \mid N=m) F(dz) \mathrm{P}(N=m)$$

$$= 2n(a_n x)^{-(\alpha+\gamma)}$$

$$\times \int_0^{a_n^\beta} z^{\alpha+\gamma} F(dz) \sum_{m=1}^{\infty} \sum_{k=1}^{m} m^{\alpha+\gamma} E(Y_k^{\alpha+\gamma} \mid N=m) \mathrm{P}(N=m)$$

$$= 2n(a_n x)^{-(\alpha+\gamma)} \int_0^{a_n^\beta} z^{\alpha+\gamma} F(dz) \mathrm{E}\left(N^{\alpha+\gamma} \sum_{k=1}^{N} Y_k^{\alpha+\gamma}\right)$$

$$\leq C n a_n^{-(\alpha+\gamma)} \int_0^{a_n^\beta} z^{\alpha+\gamma} F(dz).$$

Finally, using that $\overline{F}(z) = z^{-\alpha} L(z)$ for some slowly varying function $L$, integration by parts and the Karamata theorem, this last expression equals

$$C n a_n^{-(\alpha+\gamma)} \left(\int_0^{a_n^\beta} (\alpha+\gamma) z^{\gamma-1} L(z) \, dz - a_n^{\beta(\alpha+\gamma)} a_n^{-\beta\alpha} L(a_n^\beta)\right)$$

$$\sim C n a_n^{-(\alpha+\gamma)} \left(\frac{\alpha+\gamma}{\gamma} a_n^{\beta\gamma} L(a_n^\beta) - a_n^{\beta\gamma} L(a_n^\beta)\right)$$

$$= C \frac{\alpha}{\gamma} n a_n^{-(\alpha+\gamma(1-\beta))} L(a_n^\beta) \to 0,$$

as $n \to \infty$. Here $c_n \sim d_n$ means that $c_n/d_n \to 1$ as $n \to \infty$. In the last step we used that $(a_n)$ is regularly varying with index $1/\alpha$. This completes the proof. □

LEMMA 5.4. *Let $(Z_k)$ be an i.i.d. sequence of nonnegative random variables with $Z_1 \in \mathrm{RV}_\alpha((a_n), \mu, (0,\infty])$ and let $N$ be a $Po(\lambda)$-distributed random* variable independent of $(Z_k)$. Let $\beta \in (1/2, 1)$ and $M_n = \sum_{k=1}^{N} 1_{(a_n^\beta, \infty)}(Z_k)$. Then $\lim_{n \to \infty} n \, P(M_n \geq 2) = 0$.

PROOF. The probability generating function of $M_n$ is $g_n(t) = \exp\{\lambda p_n \times (t-1)\}$, where $p_n = P(Z_1 > a_n^\beta)$. Hence,

$$n \, P(M_n \geq 2) = n(1 - g_n(0) - g_n'(0)) = n(1 - (1 + \lambda p_n) \exp\{-\lambda p_n\})$$
$$\sim n(\lambda^2 p_n^2 / 2 + o(p_n^2)),$$

as $n \to \infty$. Since the sequence $(a_n^\beta)$ is regularly varying with index $\beta/\alpha$, for some slowly varying function $L$,

$$n p_n^2 = n(n^{-\beta} L(n))^2 = n^{1-2\beta} L^2(n) \to 0,$$

as $n \to \infty$. □

LEMMA 5.5. Let $\alpha > 0$ and let the sequence $(a_n)$ be regularly varying at infinity with index $1/\alpha$. Let $\widetilde{\mathbf{X}}$ be a Lévy process for which the Euclidean norm of each jump is bounded by 1 and let $\mathbf{Y}$ be a predictable càglàd process satisfying $E(|\mathbf{Y}|_\infty^{\alpha+\delta}) < \infty$ for some $\delta > 0$. Then $\lim_{n \to \infty} n \, P(|(\mathbf{Y} \cdot \widetilde{\mathbf{X}})|_\infty > a_n) = 0$.

PROOF. Let $\boldsymbol{\mu} = E(\widetilde{\mathbf{X}}_1)$. Then $\mathbf{M}_t = \widetilde{\mathbf{X}}_t - \boldsymbol{\mu} t$ is a martingale and

$$n \, P(a_n^{-1} |(\mathbf{Y} \cdot \widetilde{\mathbf{X}})|_\infty > \varepsilon) \leq n \, P(a_n^{-1} |\boldsymbol{\mu}| |\mathbf{Y}|_\infty > \varepsilon/2)$$
$$+ n \, P(a_n^{-1} |(\mathbf{Y} \cdot \mathbf{M})|_\infty > \varepsilon/2).$$

Let $r = \alpha + \delta/2$ so that $E(|\mathbf{Y}|_\infty^r) < \infty$. By Markov's inequality, we have, for any $\varepsilon > 0$,

$$\limsup_{n \to \infty} n \, P(a_n^{-1} |\boldsymbol{\mu}| |\mathbf{Y}|_\infty > \varepsilon/2) \leq \limsup_{n \to \infty} n a_n^{-r} (\varepsilon/2)^{-r} |\boldsymbol{\mu}|^r E(|\mathbf{Y}|_\infty^r) = 0$$

and

$$n \, P(a_n^{-1} |(\mathbf{Y} \cdot \mathbf{M})|_\infty > \varepsilon/2) \leq n a_n^{-r} (\varepsilon/2)^{-r} E(|(\mathbf{Y} \cdot \mathbf{M})|_\infty^r).$$

We will consider two cases: $\alpha \geq 1$ and $\alpha < 1$.

Assume first that $\alpha \geq 1$. Since $r > \alpha$, the claim follows if $E(|(\mathbf{Y} \cdot \mathbf{M})|_\infty^r) < \infty$. The Burkholder–Davis–Gundy inequalities (e.g., [25], page 193) and Hölder's inequality with $p = (\alpha + \delta)/(\alpha + \delta/2)$ and $q = (1 - 1/p)^{-1}$ give

$$E(|(\mathbf{Y} \cdot \mathbf{M})|_\infty^r) \leq C_r \, E\left(\left[\int_0^1 \mathbf{Y}_s^2 \, d[\mathbf{M}, \mathbf{M}]_s\right]^{r/2}\right)$$
$$\leq C_r \, E(|\mathbf{Y}|_\infty^r [\mathbf{M}, \mathbf{M}]_1^{r/2})$$
$$\leq C_r \, E(|\mathbf{Y}|_\infty^{rp})^{1/p} E([\mathbf{M}, \mathbf{M}]_1^{rq/2})^{1/q}.$$

The first factor is finite by assumption (since $rp = \alpha + \delta$) and, for some $\sigma \geq 0$,

$$[\mathbf{M}, \mathbf{M}]_t = \sigma^2 t + \sum_{0 \leq s \leq t} (\Delta \widetilde{\mathbf{X}}_s)^2, \qquad t \in [0,1],$$

which is a Lévy process with bounded jumps. Hence, by Theorem 34, page 25, in [25], $[\mathbf{M}, \mathbf{M}]_1$ has finite moments of all orders.

Assume now that $\alpha < 1$. Define the processes $\mathbf{Z}_n$ and $\widetilde{\mathbf{Z}}_n$ by

$$\mathbf{Z}_n(s) = \mathbf{Y}_s 1_{(a_n, \infty)}\left(\sup_{u \in [0,s]} |\mathbf{Y}_u|\right), \qquad s \in [0,1],$$

$$\widetilde{\mathbf{Z}}_n(s) = \mathbf{Y}_s 1_{[0, a_n]}\left(\sup_{u \in [0,s]} |\mathbf{Y}_u|\right), \qquad s \in [0,1],$$

and note that $\mathbf{Y} = \mathbf{Z}_n + \widetilde{\mathbf{Z}}_n$ so that $(\mathbf{Y} \cdot \mathbf{M}) = (\mathbf{Z}_n \cdot \mathbf{M}) + (\widetilde{\mathbf{Z}}_n \cdot \mathbf{M})$. Moreover,

$$n\,\mathrm{P}(a_n^{-1}|(\mathbf{Y} \cdot \mathbf{M})|_\infty > \varepsilon/2)$$
$$\leq n\,\mathrm{P}(a_n^{-1}|(\mathbf{Z}_n \cdot \mathbf{M})|_\infty > \varepsilon/4) + n\,\mathrm{P}(a_n^{-1}|(\widetilde{\mathbf{Z}}_n \cdot \mathbf{M})|_\infty > \varepsilon/4)$$
$$\leq n\,\mathrm{P}(|\mathbf{Y}|_\infty > a_n) + n\,\mathrm{P}(a_n^{-1}|(\widetilde{\mathbf{Z}}_n \cdot \mathbf{M})|_\infty > \varepsilon/4).$$

Markov's inequality yields $\lim_{n \to \infty} n\,\mathrm{P}(|\mathbf{Y}|_\infty > a_n) = 0$. The Burkholder–Davis–Gundy and Hölder inequalities yield

$$\mathrm{E}(|(\widetilde{\mathbf{Z}}_n \cdot \mathbf{M})|_\infty^2) \leq C_2\,\mathrm{E}(|\widetilde{\mathbf{Z}}_n|_\infty^{2p})^{1/p}\,\mathrm{E}([\mathbf{M}, \mathbf{M}]^{2q/2})^{1/q} = K\,\mathrm{E}(|\widetilde{\mathbf{Z}}_n|_\infty^{2p})^{1/p},$$

for any $p > 1$, $q = (1 - 1/p)^{-1}$, where $K \in (0, \infty)$ is a constant. If we put $r = \alpha + \delta/2$ and $p = (\alpha + \delta)/(\alpha + \delta/2)$, then we obtain

$$n\,\mathrm{P}(a_n^{-1}|(\widetilde{\mathbf{Z}}_n \cdot \mathbf{M})|_\infty > \varepsilon/4) \leq n a_n^{-r} a_n^{-2+r}(\varepsilon/4)^{-2}\,\mathrm{E}(|(\widetilde{\mathbf{Z}}_n \cdot \mathbf{M})|_\infty^2)$$
$$\leq K(\varepsilon/4)^{-2} n a_n^{-r}\,\mathrm{E}(a_n^{-p(2-r)}|\widetilde{\mathbf{Z}}_n|_\infty^{2p})^{1/p}.$$

Note that $\lim_{n \to \infty} n a_n^{-r} = 0$. For the expectation above, we have, with $F$ denoting the distribution function of $|\mathbf{Y}|_\infty$,

$$\begin{aligned}
\mathrm{E}(a_n^{-p(2-r)}|\widetilde{\mathbf{Z}}_n|_\infty^{2p}) &= \int_0^{a_n} x^{rp}(x/a_n)^{2p-rp}\,dF(x) \\
&\leq \int_0^{a_n} x^{rp}\,dF(x) \\
&= \mathrm{E}(|\mathbf{Y}|_\infty^{rp} 1_{[0, a_n]}(|\mathbf{Y}|_\infty)) \\
&\to \mathrm{E}(|\mathbf{Y}|_\infty^{rp}) < \infty.
\end{aligned}$$

The conclusion follows. $\square$

PROPOSITION 5.1. *Assume the hypotheses of Theorem 3.4 and let* **J** *be the compound Poisson part as in* (1.4). *Then,* $(\mathbf{Y} \cdot \mathbf{J}) \in \mathrm{RV}((a_n), m^*, \overline{\mathbb{D}}_0)$, *where*

$$m^*(B) = \mathrm{E}(\mu\{\mathbf{x} \in \overline{\mathbb{R}}_0^d : \mathbf{Y}_V \mathbf{x} 1_{[V,1]} \in B\}),$$

*where $V$ is uniformly distributed on $[0,1)$ and independent of* **Y**.

PROOF. Take constants $\beta \in (1/2, 1)$ and $C > 0$ (we will eventually let $C \to \infty$). Put, as in the proof of Theorem 3.4,

$$\mathbf{J}_n = \sum_{k=1}^{N_1} \mathbf{Z}_k 1_{(a_n^\beta, \infty)}(|\mathbf{Z}_k|) 1_{[\tau_k, 1]} \quad \text{and} \quad M_n = \sum_{k=1}^{N_1} 1_{(a_n^\beta, \infty)}(|\mathbf{Z}_k|).$$

With this notation, we may also write

$$\mathbf{J}_n = \sum_{k=1}^{M_n} \mathbf{Z}_k^{(n)} 1_{[\tau_k^{(n)}, 1]},$$

where $\mathbf{Z}_k^{(n)}$ is the $k$th jump with norm larger than $a_n^\beta$, and $\tau_k^{(n)}$ is the time of that jump. Take closed $B \in \mathcal{B}(\overline{\mathbb{D}}_0) \cap \mathbb{D}$ bounded away from **0** and let $d_B = \inf\{|\mathbf{x}|_\infty : \mathbf{x} \in B\}$. For $\varepsilon \in (0, d_B)$, let $B_\varepsilon = \{\mathbf{x} \in \mathbb{D} : d^\circ(\mathbf{x}, B) \le \varepsilon\}$. Note that $B_\varepsilon$ is closed. By the Portmanteau theorem and Remark 4.1 it is sufficient to prove that

(5.3) $$\limsup_{n \to \infty} n \, \mathrm{P}(a_n^{-1}(\mathbf{Y} \cdot \mathbf{J}) \in B) \le m^*(B)$$

and for $\delta > 0$,

(5.4) $$\lim_{n \to \infty} n \, \mathrm{P}(|(\mathbf{Y} \cdot \mathbf{J})|_\infty \ge a_n \delta) = m^*(B_{0,\delta}^c).$$

The outline of the proof of (5.3) is as follows:

(i) First we will show that

$$\limsup_{n \to \infty} n \, \mathrm{P}(a_n^{-1}(\mathbf{Y} \cdot \mathbf{J}) \in B)$$

$$\le \lim_{C \to \infty} \limsup_{n \to \infty} n \, \mathrm{P}(\mathbf{Y}_{\tau_1^{(n)}} \mathbf{Z}_1^{(n)} 1_{[\tau_1^{(n)}, 1]} \in a_n B_\varepsilon, M_n = 1, |\mathbf{Y}_{\tau_1^{(n)}}| \le C).$$

(ii) Then we will show that

$$\limsup_{n \to \infty} n \, \mathrm{P}(\mathbf{Y}_{\tau_1^{(n)}} \mathbf{Z}_1^{(n)} 1_{[\tau_1^{(n)}, 1]} \in a_n B_\varepsilon, M_n = 1, |\mathbf{Y}_{\tau_1^{(n)}}| \le C)$$

$$\le \int_{B_{\mathbf{0},C} \times [0,1]} \mu\{\mathbf{x} \in \mathbb{R}^d : \mathbf{xy} 1_{[t,1]} \in B_\varepsilon\} \rho(d(\mathbf{y}, t)),$$

where $\rho(A \times T) = \mathrm{P}((\mathbf{Y}_V, V) \in A \times T)$ with $V$ uniformly distributed on $[0,1)$ and independent of **Y**.

Finally, letting $C \to \infty$ and then $\varepsilon \downarrow 0$, the conclusion follows.

Let us first prove (i). We have,

$$\{a_n^{-1}(\mathbf{Y} \cdot \mathbf{J}) \in B\} = \{a_n^{-1}(\mathbf{Y} \cdot \mathbf{J}) \in B, |(\mathbf{Y} \cdot (\mathbf{J} - \mathbf{J}_n))|_\infty > a_n \varepsilon\}$$
$$\cup \{a_n^{-1}(\mathbf{Y} \cdot \mathbf{J}) \in B, |(\mathbf{Y} \cdot (\mathbf{J} - \mathbf{J}_n))|_\infty \le a_n \varepsilon\}$$
$$\subset \{|(\mathbf{Y} \cdot (\mathbf{J} - \mathbf{J}_n))|_\infty > a_n \varepsilon\} \cup \{(\mathbf{Y} \cdot \mathbf{J}_n) \in a_n B_\varepsilon\}.$$

By Lemma 5.3 and Remark 5.1, $\lim_{n \to \infty} n \operatorname{P}(|(\mathbf{Y} \cdot (\mathbf{J} - \mathbf{J}_n))|_\infty > a_n \varepsilon) = 0$. For the second term, we have

$$(5.5) \qquad n \operatorname{P}((\mathbf{Y} \cdot \mathbf{J}_n) \in a_n B_\varepsilon) = n \operatorname{P}((\mathbf{Y} \cdot \mathbf{J}_n) \in a_n B_\varepsilon, M_n = 1)$$
$$(5.6) \qquad\qquad\qquad\qquad + n \operatorname{P}((\mathbf{Y} \cdot \mathbf{J}_n) \in a_n B_\varepsilon, M_n \ge 2).$$

By Lemma 5.4, the term (5.6) converges to zero as $n \to \infty$. It remains to consider (5.5). We have

$$n \operatorname{P}((\mathbf{Y} \cdot \mathbf{J}_n) \in a_n B_\varepsilon, M_n = 1)$$
$$(5.7) \qquad = n \operatorname{P}(\mathbf{Y}_{\tau_1^{(n)}} \mathbf{Z}_1^{(n)} 1_{[\tau_1^{(n)}, 1]} \in a_n B_\varepsilon, M_n = 1, |\mathbf{Y}_{\tau_1^{(n)}}| > C)$$
$$(5.8) \qquad + n \operatorname{P}(\mathbf{Y}_{\tau_1^{(n)}} \mathbf{Z}_1^{(n)} 1_{[\tau_1^{(n)}, 1]} \in a_n B_\varepsilon, M_n = 1, |\mathbf{Y}_{\tau_1^{(n)}}| \le C).$$

Since $\mathbf{Y}_{\tau_1^{(n)}}$ and $\mathbf{Z}_1^{(n)}$ are independent, we apply Lemma 5.2 to (5.7) and obtain

$$n \operatorname{P}(\mathbf{Y}_{\tau_1^{(n)}} \mathbf{Z}_1^{(n)} 1_{[\tau_1^{(n)}, 1]} \in a_n B_\varepsilon, M_n = 1, |\mathbf{Y}_{\tau_1^{(n)}}| > C)$$
$$(5.9) \qquad \le n \operatorname{P}(|\mathbf{Y}_{\tau_1^{(n)}}| 1_{(C, \infty)}(|\mathbf{Y}_{\tau_1^{(n)}}|) |\mathbf{Z}_1^{(n)}| > a_n(d_B - \varepsilon))$$
$$\le 2n \operatorname{P}(|\mathbf{Y}_{\tau_1^{(n)}}| 1_{(C, \infty)}(|\mathbf{Y}_{\tau_1^{(n)}}|) |\widetilde{\mathbf{Z}}| 1_{(a_n^\beta, \infty)}(|\widetilde{\mathbf{Z}}|) > a_n(d_B - \varepsilon)),$$

where $\widetilde{\mathbf{Z}} \stackrel{d}{=} \mathbf{Z}_1$ and is independent of $\mathbf{Y}$ and $\mathbf{J}$. Hence, Breiman's result (3.1) can be applied to show that (5.9) satisfies

$$2n \operatorname{P}(|\mathbf{Y}_{\tau_1^{(n)}}| 1_{(C, \infty)}(|\mathbf{Y}_{\tau_1^{(n)}}|) |\widetilde{\mathbf{Z}}| 1_{(a_n^\beta, \infty)}(|\widetilde{\mathbf{Z}}|) > a_n(d_B - \varepsilon))$$
$$\le 2n \operatorname{P}(|\mathbf{Y}|_\infty 1_{(C, \infty)}(|\mathbf{Y}|_\infty) |\widetilde{\mathbf{Z}}| > a_n(d_B - \varepsilon))$$
$$\to 2 \operatorname{E}(|\mathbf{Y}|_\infty^\alpha 1_{(C, \infty)}(|\mathbf{Y}|_\infty)) \mu(B_{\mathbf{0}, d_B - \varepsilon}^c).$$

Finally, letting $C \to \infty$, the last expression converges to 0. This completes the proof of (i).

(ii) We now study (5.8). Set $\Gamma(C) = \{(\mathbf{y}, t) \in \mathbb{R}^d \times [0, 1) : |\mathbf{y}| \le C\}$. Conditioning on $(\mathbf{Y}_{\tau_1^{(n)}}, \tau_1^{(n)})$, we get

$$n \operatorname{P}(\mathbf{Y}_{\tau_1^{(n)}} \mathbf{Z}_1^{(n)} 1_{[\tau_1^{(n)}, 1]} \in a_n B_\varepsilon, M_n = 1, |\mathbf{Y}_{\tau_1^{(n)}}| \le C)$$

$$= \int_{\Gamma(C)} n\, \mathrm{P}(\mathbf{y}\mathbf{Z}_1^{(n)} 1_{[t,1]} \in a_n B_\varepsilon \mid M_n = 1)$$
$$\times \mathrm{P}((\mathbf{Y}_{\tau_1^{(n)}}, \tau_1^{(n)}) \in d(\mathbf{y}, t), M_n = 1)$$
$$= \int_{\Gamma(C)} \underbrace{n\, \mathrm{P}(\mathbf{y}\mathbf{Z}_1^{(n)} 1_{[t,1]} \in a_n B_\varepsilon, M_n = 1)}_{f_n(\mathbf{y},t)}$$
$$\times \underbrace{\frac{\mathrm{P}((\mathbf{Y}_{\tau_1^{(n)}}, \tau_1^{(n)}) \in d(\mathbf{y}, t), M_n = 1)}{\mathrm{P}(M_n = 1)}}_{\rho_n(d(\mathbf{y},t))}.$$

For $\mathbf{w} \in \mathbb{D}$, we denote by $\varphi_\mathbf{w} : \mathbb{R}^d \to \mathbb{D}$ the function given by $\varphi_\mathbf{w}(\mathbf{x}) = \mathbf{x}\mathbf{w}$. By Theorem 4.2 in [30], multiplication $\psi : \mathbb{D} \times \mathbb{D} \to \mathbb{D}$ given by

$$\psi(\mathbf{w}, \mathbf{z})_t = \mathbf{w}_t \mathbf{z}_t = (w_t^{(1)} z_t^{(1)}, \ldots, w_t^{(d)} z_t^{(d)}), \qquad t \in [0,1],$$

is continuous at those $(\mathbf{w}, \mathbf{z}) \in \mathbb{D} \times \mathbb{D}$ for which $\mathrm{Disc}(\mathbf{w}) \cap \mathrm{Disc}(\mathbf{z}) = \varnothing$. Moreover, $h : \mathbb{R}^d \to \mathbb{D}$ given by $h(\mathbf{x}) = \mathbf{x} 1_{[0,1]}$ is continuous and $\mathrm{Disc}(h(\mathbf{x})) = \varnothing$ for every $\mathbf{x} \in \mathbb{R}^d$. Hence, $\varphi_\mathbf{w}(\cdot) = \psi(h(\cdot), \mathbf{w})$ is continuous for every $\mathbf{w} \in \mathbb{D}$. We will show the following:

(a) $\limsup_{n \to \infty} \sup_{(\mathbf{y},t) \in \Gamma(C)} (f_n(\mathbf{y}, t) - f(\mathbf{y}, t)) \leq 0$, $f(\mathbf{y}, t) = \mu \circ \varphi_{\mathbf{y}1_{[t,1]}}^{-1}(B_\varepsilon)$.

(b) $\rho_n \xrightarrow{w} \rho$, $\rho(A \times T) = \mathrm{P}((\mathbf{Y}_V, V) \in A \times T)$, $V$ uniformly distributed on $[0, 1)$ and independent of $\mathbf{Y}$.

(c) $\limsup_{n \to \infty} \int_{\Gamma(C)} f(\mathbf{y}, t) \rho_n(d(\mathbf{y}, t)) \leq \int_{\Gamma(C)} f(\mathbf{y}, t) \rho(d(\mathbf{y}, t))$.

Using (a)–(c), it follows that

$$\limsup_{n \to \infty} \left( \int_{\Gamma(C)} f_n(\mathbf{y}, t) \rho_n(d\mathbf{y} \times dt) - \int_{\Gamma(C)} f(\mathbf{y}, t) \rho(d(\mathbf{y}, t)) \right)$$
$$\leq \limsup_{n \to \infty} \sup_{(\mathbf{y},t) \in \Gamma(C)} (f_n(\mathbf{y}, t) - f(\mathbf{y}, t)) \rho_n(\Gamma(C))$$
$$+ \limsup_{n \to \infty} \int_{\Gamma(C)} f(\mathbf{y}, t) (\rho_n(d(\mathbf{y}, t)) - \rho(d(\mathbf{y}, t)))$$
$$\leq 0.$$

This proves (ii). We start by showing (a). Since $B_\varepsilon$ is closed, it follows, by continuity, that $\varphi_{\mathbf{y}1_{[t,1]}}^{-1}(B_\varepsilon)$ is closed. For large enough $n$, we have

$$f_n(\mathbf{y}, t) = n\, \mathrm{P}(\mathbf{y}\mathbf{Z}_1^{(n)} 1_{[t,1]} \in a_n B_\varepsilon, M_n = 1)$$
$$\leq n\, \mathrm{P}(\mathbf{y}\mathbf{Z}_1 1_{[t,1]} \in a_n B_\varepsilon)$$
$$= n\, \mathrm{P}(a_n^{-1} \mathbf{Z}_1 \in \varphi_{\mathbf{y}1_{[t,1]}}^{-1}(B_\varepsilon)).$$

Hence, by the continuous mapping theorem and the Portmanteau theorem,

$$\limsup_{n\to\infty} f_n(\mathbf{y}, t) \leq f(\mathbf{y}, t).$$

Recall the uniformity of regular variation: if $\delta > 0$ and the distribution $F$ on $\mathbb{R}^d$ is regularly varying, that is, $F \in \mathrm{RV}((a_n), \mu, \overline{\mathbb{R}}_0^d)$, then, for each $\eta > 0$, there exists $N(\eta)$ such that, for $n \geq N(\eta)$ and each closed set $B \subset B_{\mathbf{0},\delta}^c$, we have

$$nF(a_n B) \leq \mu(B) + \eta.$$

In our setting we have for each $(\mathbf{y}, t)$ with $|\mathbf{y}| \leq C$ and $t \in [0,1)$ that $\{\mathbf{z} \in \mathbb{R}^d : \mathbf{y}\mathbf{z}1_{[t,1]} \in B_\varepsilon\} \subset \{\mathbf{z} \in \mathbb{R}^d : |\mathbf{z}| > (d_B - \varepsilon)/C\}$. Hence, for each $\eta > 0$, there exists $N(\eta)$ such that, for $n \geq N(\eta)$

$$n\mathrm{P}(\mathbf{y}\mathbf{Z}_1 1_{[t,1]} \in a_n B_\varepsilon) \leq \mu \circ \varphi_{\mathbf{y}1_{[t,1]}}^{-1}(B_\varepsilon) + \eta,$$

uniformly on $\Gamma(C)$. That is,

$$\limsup_{n\to\infty} \sup_{(\mathbf{y},t)\in\Gamma(C)} (f_n(\mathbf{y}, t) - f(\mathbf{y}, t)) \leq 0.$$

For (b), we have the following. Let $A \times T \in \mathcal{B}(\mathbb{R}^d \times [0,1))$ be a $\rho$-continuity set. Conditioning on $\tau_1^{(n)}$ and using that $\mathbf{Y}$ is predictable, we have

$$\mathrm{P}((\mathbf{Y}_{\tau_1^{(n)}}, \tau_1^{(n)}) \in A \times T \mid M_n = 1)$$

$$= \int_T \mathrm{P}(\mathbf{Y}_t \in A \mid \tau^{(n)} = t, M_n = 1)\, dt$$

$$= \int_T \mathrm{P}\left(\mathbf{Y}_t \in A \mid \sup_{s<t} |\Delta \mathbf{X}_s| \leq a_n^\beta\right) dt.$$

Since $\lim_{n\to\infty} \mathrm{P}(\mathbf{Y}_t \in A \mid \sup_{s<t} |\Delta \mathbf{X}_s| \leq a_n^\beta) = \mathrm{P}(\mathbf{Y}_t \in A)$ for all but at most countably many $t \in [0,1)$, the dominated convergence theorem yields

$$\lim_{n\to\infty} \mathrm{P}((\mathbf{Y}_{\tau_1^{(n)}}, \tau_1^{(n)}) \in A \times T \mid M_n = 1) = \int_T \mathrm{P}(\mathbf{Y}_t \in A)\, dt$$

$$= \mathrm{P}((\mathbf{Y}_V, V) \in A \times T),$$

where $V$ is uniformly distributed on $[0,1)$ and independent of $\mathbf{Y}$. This proves (b).

Finally, we have to prove (c). Given $\mathbf{x} \in \mathbb{R}^d$, denote by $\widetilde{\varphi}_{\mathbf{x}} : \mathbb{R}^d \times [0,1) \to \mathbb{D}$ the mapping $\widetilde{\varphi}_{\mathbf{x}}(\mathbf{y}, t) = \mathbf{y}\mathbf{x}1_{[t,1]}$. For each $\mathbf{x} \in \mathbb{R}^d$, the mapping $\widetilde{\varphi}_{\mathbf{x}}$ is continuous. Let $(\mathbf{U}_n, V_n)$ and $(\mathbf{U}, V)$ be random vectors with distribution $\rho_n$ and $\rho$, respectively. We have

$$\int_{\Gamma(C)} f(\mathbf{y}, t) \rho_n(d(\mathbf{y}, t))$$

$$= \mathrm{E}(\mu \circ \varphi^{-1}_{\mathbf{U}_n 1_{[V_n,1]}}(B_\varepsilon) 1_{\Gamma(C)}(\mathbf{U}_n, V_n))$$

$$= \int_\Omega \int_{\mathbb{R}^d \setminus \{\mathbf{0}\}} 1_{B_\varepsilon}(\mathbf{x}\mathbf{U}_n(\omega) 1_{[V_n(\omega),1]}) 1_{\Gamma(C)}(\mathbf{U}_n(\omega), V_n(\omega)) \mu(d\mathbf{x}) \mathrm{P}(d\omega)$$

$$= \int_{\mathbb{R}^d \setminus \{\mathbf{0}\}} \mathrm{E}(1_{B_\varepsilon}(\mathbf{x}\mathbf{U}_n 1_{[V_n,1]}) 1_{\Gamma(C)}(\mathbf{U}_n, V_n)) \mu(d\mathbf{x})$$

$$= \int_{\mathbb{R}^d \setminus \{\mathbf{0}\}} \mathrm{P}(\mathbf{x}\mathbf{U}_n 1_{[V_n,1]} \in B_\varepsilon, (\mathbf{U}_n, V_n) \in \Gamma(C)) \mu(d\mathbf{x})$$

$$= \int_{\mathbb{R}^d \setminus \{\mathbf{0}\}} \mathrm{P}((\mathbf{U}_n, V_n) \in \widetilde{\varphi}_{\mathbf{x}}^{-1}(B_\varepsilon) \cap \Gamma(C)) \mu(d\mathbf{x}).$$

By (b), $(\mathbf{U}_n, V_n) \xrightarrow{d} (\mathbf{U}, V)$. Moreover, since $B_\varepsilon$ is closed, it follows that $\widetilde{\varphi}_{\mathbf{x}}^{-1}(B_\varepsilon)$ is closed. Since $\Gamma(C)$ is closed, also $\widetilde{\varphi}_{\mathbf{x}}^{-1}(B_\varepsilon) \cap \Gamma(C)$ is closed. Hence, by the Portmanteau theorem,

$$\limsup_{n \to \infty} \mathrm{P}((\mathbf{U}_n, V_n) \in \widetilde{\varphi}_{\mathbf{x}}^{-1}(B_\varepsilon) \cap \Gamma(C)) \leq \mathrm{P}((\mathbf{U}, V) \in \widetilde{\varphi}_{\mathbf{x}}^{-1}(B_\varepsilon) \cap \Gamma(C)),$$

and we arrive at

$$\limsup_{n \to \infty} \int_{\Gamma(C)} f(\mathbf{y}, t) \rho_n(d(\mathbf{y}, t))$$

$$\leq \limsup_{n \to \infty} \int_{\mathbb{R}^d \setminus \{\mathbf{0}\}} \mathrm{P}((\mathbf{U}_n, V_n) \in \widetilde{\varphi}_{\mathbf{x}}^{-1}(B_\varepsilon) \cap \Gamma(C)) \mu(d\mathbf{x})$$

$$\leq \int_{\mathbb{R}^d \setminus \{\mathbf{0}\}} \limsup_{n \to \infty} \mathrm{P}((\mathbf{U}_n, V_n) \in \widetilde{\varphi}_{\mathbf{x}}^{-1}(B_\varepsilon) \cap \Gamma(C)) \mu(d\mathbf{x})$$

$$\leq \int_{\mathbb{R}^d \setminus \{\mathbf{0}\}} \mathrm{P}((\mathbf{U}, V) \in \widetilde{\varphi}_{\mathbf{x}}^{-1}(B_\varepsilon) \cap \Gamma(C)) \mu(d\mathbf{x})$$

$$= \int_{\Gamma(C)} f(\mathbf{y}, t) \rho(d(\mathbf{y}, t)).$$

The interchange of the limit and the integral is allowed if there is a function $g$ such that $\mathrm{P}((\mathbf{U}_n, V_n) \in \widetilde{\varphi}_{\mathbf{x}}^{-1}(B_\varepsilon) \cap \Gamma(C)) \leq g(\mathbf{x})$ and $\int_{\mathbb{R}^d \setminus \{\mathbf{0}\}} g(\mathbf{x}) \mu(d\mathbf{x}) < \infty$. We have

$$\mathrm{P}((\mathbf{U}_n, V_n) \in \widetilde{\varphi}_{\mathbf{x}}^{-1}(B_\varepsilon) \cap \Gamma(C)) \leq \mathrm{P}(|\mathbf{x}\mathbf{U}_n| > d_B - \varepsilon, |\mathbf{U}_n| \leq C)$$

$$\leq \mathrm{P}(|\mathbf{U}_n| \in ((d_B - \varepsilon)/|\mathbf{x}|, C])$$

$$\leq 1_{((d_B - \varepsilon)/C, \infty)}(|\mathbf{x}|)$$

and $\mu\{\mathbf{x} \in \mathbb{R}^d : |\mathbf{x}| > (d_B - \varepsilon)/C\} < \infty$. Hence, we may take $g(\mathbf{x}) = 1_{((d_B - \varepsilon)/C, \infty)}(|\mathbf{x}|)$. This concludes the proof of (c) and hence of (ii) and the proof of (5.3) is complete. It remains to prove (5.4). From (5.3) we have

the upper bound

$$\limsup_{n\to\infty} n\, P(|(\mathbf{Y}\cdot\mathbf{J})|_\infty \geq a_n\delta) \leq m^*(B_{0,\delta}^c).$$

The proof of the lower bound

$$\liminf_{n\to\infty} n\, P(|(\mathbf{Y}\cdot\mathbf{J})|_\infty > a_n\delta) \geq m^*(B_{0,\delta}^c)$$

is similar, replacing lim sup by lim inf, replacing closed sets by open sets, changing the direction of inequalities and using that $m^*(\partial B_{0,\delta}^c) = 0$. This completes the proof. $\square$

**Acknowledgment.** The authors want to thank the referee for comments leading to Remark 3.2.

DIVISION OF APPLIED MATHEMATICS
BROWN UNIVERSITY
182 GEORGE ST
PROVIDENCE, RHODE ISLAND 02912
USA
E-MAIL: henrik_hult@brown.edu

DEPARTMENT OF MATHEMATICS
KTH
SE-100 44 STOCKHOLM
SWEDEN
E-MAIL: lindskog@math.kth.se